\documentclass[a4paper,11pt]{amsart}
\usepackage{amsmath}
\usepackage{amssymb}
\usepackage{amsthm}
\usepackage[all]{xy}
\usepackage[latin1]{inputenc}        % accents
\usepackage[dvips]{graphics}
\usepackage[dvips]{graphicx}
\usepackage{amscd}
\usepackage{mathrsfs}
\usepackage[mathcal]{eucal}
\usepackage{float}
\usepackage{fullpage}
\author{Francesco Polizzi}
\address{Dipartimento di Matematica, Università della Calabria, Via Pietro Bucci,
87036 Arcavacata di Rende (CS), Italy.}
\email{polizzi@mat.unical.it}
\title[Standard isotrivial fibrations]{Standard isotrivial fibrations with $p_g=q=1$}
\date{\today}

\newtheorem{inizio}{Lemma}[section]
\newtheorem{theorem}[inizio]{Theorem}
\newtheorem{corollary}[inizio]{Corollary}
\newtheorem{proposition}[inizio]{Proposition}
\newtheorem{lemma}[inizio]{Lemma}

\newtheorem{definition}[inizio]{Definition}
\newtheorem*{teo}{Main Theorem}

\theoremstyle{definition}
\newtheorem{remark}[inizio]{Remark}

\newcommand{\lr}{\longrightarrow}

\newcommand{\mO}{\mathcal{O}}
\newcommand{\mZ}{\mathbb{Z}}

\newcommand{\Sn}{\sum_{j=1}^s \left(1- \frac{1}{ \;n_j} \right)}

\begin{document}

\subjclass{14J29 (primary), 14Q99, 20F05} \keywords{Surfaces of general
type, isotrivial fibrations, actions of finite groups}

\maketitle

\abstract A smooth, projective surface $S$ of general type is said to be a
\emph{standard isotrivial fibration} if there exist a finite group
$G$ acting faithfully on two smooth projective curves $C$ and
$F$ so that $S$ is isomorphic to the minimal desingularization of
$T:=(C \times F)/G$. If $T$ is smooth then $S=T$ is called a
$\emph{quasi-bundle}$. In this paper we classify the standard
isotrivial fibrations with $p_g=q=1$ which are not quasi-bundles,
assuming that all the singularities of $T$ are rational double
points. As a by-product, we provide several new examples of minimal
surfaces of general type with $p_g=q=1$ and $K_S^2=4,6$.
\endabstract

%\tableofcontents
\section{Introduction}
Recently, there has been considerable interest in understanding the
geometry of complex projective surfaces with small birational
invariants, and in particular of surfaces with $p_g=q=1$. Any
surface $S$ of general type verifies $\chi(\mO_S) >0$, hence
$q(S)>0$ implies $p_g(S)>0$. It follows that the surfaces of general
type with $p_g=q=1$ are the irregular ones
 with the lowest geometric genus, hence it would be important to
 achieve their complete classification. So far, this has been
 obtained only in the cases $K_S^2=2,3$ (see \cite{Ca},
 \cite{CaCi91}, \cite{CaCi93}, \cite{Pol05}, \cite{CaPi06}).
If $S$ is any surface with $q=1$, its Albanese map $\alpha \colon S
\lr E$ is a fibration over an elliptic curve $E$; we denote by
$g_{\textrm{alb}}$ the genus of the general fibre of $\alpha$. The
universal property of the Albanese morphism implies that $\alpha$ is
the unique fibration on $S$ with irrational base. As the title
suggests, this paper considers surfaces with $p_g=q=1$ which are
\emph{standard isotrivial fibrations}. This means that there exists
a finite group $G$ which acts faithfully on two smooth projective
curves $C$ and $F$ so that $S$ is isomorphic to the minimal
desingularization of $T:=(C \times F)/G$. If $T$ is smooth then
$S=T$ is called a $\emph{quasi-bundle}$ or a surface \emph{isogenous
to an unmixed product} (see \cite{Se90}, \cite{Se96}, \cite{Ca00}).
Quasi-bundles of general type with $p_g=q=1$ are classified in
\cite{Pol08} and \cite{CarPol}. In the present work we consider the
case where all the singularities of $T$ are rational double points
(RDPs). Our
 classification procedure combines ideas from \cite{Pol08}
 and combinatorial methods
 of finite group theory. Let $\lambda \colon S \lr T=(C \times F)/G$ be a standard
isotrivial fibration; then the two projections $\pi_C \colon C
\times F \lr C$, $\pi_F \colon C \times F \lr F$ induce
 two morphisms $\alpha \colon S \lr C/G$, $\beta \colon S \lr F/G$, whose
smooth fibres are isomorphic to $F$ and $C$, respectively. We have
$q(S)=g(C/G)+g(F/G)$, then if $q(S)=1$ we may assume that $E:=C/G$
is an elliptic curve and $F/G \cong \mathbb{P}^1$. Consequently, the
morphism $\alpha$ is the Albanese fibration of $S$ and
$g_{\textrm{alb}}=g(F)$. If $p_g(S)=q(S)=1$ and $T$
contains only RDPs, we show that $S$ is a minimal surface (Proposition
\ref{minimal}) and that $2 \leq g(F) \leq 4$. Therefore we
can use the classification of finite groups acting on Riemann
surfaces of low genus given in \cite{Br90},
\cite{KuKi90}, \cite{KuKu90}, \cite{Bre00}, \cite{Vin00}, \cite{Ki03}.
In particular we obtain
$|G| \leq 168$ and so the problem can be attacked with the computer
algebra program GAP4, whose database includes all groups of order
less than $2000$, with the exception of $1024$ (see \cite{GAP4}).
Computer algebra is
 a powerful tool when dealing with this kind of problems;
a recent example is the paper \cite{BaCaGr06}, where the MAGMA
database of finite groups (identical to the GAP4 database) is
exploited in order to achieve the classification of surfaces with
$p_g=q=0$ isogenous to a product. In our case we have tried to
minimize the amount of computer calculations, doing everything ``by
hand'' whenever possible and using GAP4 only when working with
groups of big order or cumbersome presentation. Nevertheless, the
computer's aid has been extremely useful in order to obtain some of
the non-generation results of Section \ref{non-generation} and some
of the existence results of Section \ref{general case}. The aim of
this paper is to prove the following
\begin{teo} Let $\lambda \colon S \lr T=(C \times F)/G$ be
a standard isotrivial fibration of general type with $p_g=q=1$,
which is not a quasi-bundle, and assume that $T$ contains only
$\emph{RDPs}$. Then $S$ is a minimal surface,
$K_S^2$ is even and the singularities of $T$ are
exactly $8-K_S^2$ nodes. Moreover, the occurrences for $K_S^2$, $g(F)$,
$g(C)$ and $G$ are precisely those listed in the table below.
\end{teo}
\begin{table}[ht!]
\begin{center}
\begin{tabular}{|c|c|c|c|c|}
\hline
$ $ & $ $ &  $ $ & $ $ &  \verb|IdSmall| \\
$K_S^2$  & $g(F)=g_{\textrm{alb}}$ & $g(C)$ & $G$ & \verb|Group|$(G)$  \\
\hline \hline
$6$ & $3$ & $10$ & $\textrm{SL}_2(\mathbb{F}_3)$ & $G(24,3)$ \\
\hline
$6$ & $3$ & $13$ & $\mZ_2 \ltimes (\mZ_2 \times \mZ_8)$ & $G(32,9)$ \\
\hline
$6$ & $3$ & $13$ & $\mZ_2 \ltimes D_{2,8,5}$ & $G(32,11)$ \\
\hline
$6$ & $3$ & $19$ & $G(48,33)$ & $G(48,33)$ \\
\hline
$6$ & $3$ & $19$ & $\mZ_3 \ltimes (\mZ_4)^2$ & $G(48,3)$ \\
\hline
$6$ & $3$ & $64$ & $\textrm{PSL}_2(\mathbb{F}_7)$ & $G(168,42)$ \\
\hline
$6$ & $4$ & $3$ & $D_4$ & $G(8,3)$ \\
\hline
$6$ & $4$ & $4$ & $A_4$ & $G(12,3)$ \\
\hline
$6$ & $4$ & $7$ & $D_{2,12,7}$ & $G(24,10)$ \\
\hline
$6$ & $4$ & $10$ & $\mZ_3 \times A_4$ & $G(36,11)$ \\
\hline
$6$ & $4$ & $19$ & $D_4 \ltimes (\mZ_3)^2$ & $G(72,40)$ \\
\hline
$6$ & $4$ & $31$ & $S_5$ & $G(120,34)$ \\
\hline
$4$ & $2$ & $3$ &  $\mZ_2 \times \mZ_2$ & $G(4,2)$ \\
\hline
$4$ & $2$ & $4$ & $\mZ_6$ & $G(6,2)$ \\
\hline
$4$ & $2$ & $4$ & $S_3$ & $G(6,1)$ \\
\hline
$4$ & $2$ & $5$ & $D_4$ & $G(8,3)$ \\
\hline
$4$ & $2$ & $7$ & $\mZ_2 \times \mZ_6$ & $G(12,5)$ \\
\hline
$4$ & $2$ & $7$ & $D_6$ & $G(12,4)$ \\
\hline
$4$ & $2$ & $9$ & $D_{2,8,3}$ & $G(16,8)$ \\
\hline
$4$ & $2$ & $13$ & $\mZ_2 \ltimes ((\mZ_2)^2 \times \mZ_3)$ & $G(24,8)$ \\
\hline
$4$ & $2$ & $25$ & $\textrm{GL}_2(\mathbb{F}_3)$ & $G(48,29)$ \\
\hline
$4$ & $3$ & $3$ & $D_4$ & $G(8,3)$\\
\hline
$4$ & $3$ & $4$ & $A_4$ & $G(12,3)$ \\
\hline
$4$ & $3$ & $5$ & $D_{2,8,5}$ & $G(16,6)$ \\
\hline
$4$ & $3$ & $5$ & $D_{4,4,-1}$ & $G(16,4)$ \\
\hline
$4$ & $3$ & $7$ & $\mZ_2 \times A_4$ & $G(24,13)$\\
\hline
$2$ & $2$ & $3$ & $Q_8$ & $G(8,4)$ \\
\hline
$2$ & $2$ & $3$ & $D_4$ & $G(8,3)$ \\
\hline
\end{tabular}
\end{center}
\end{table}
Here \verb|IdSmallGroup|$(G)$ denotes the label of $G$ in the GAP4
database of small groups. For instance,
\verb|IdSmallGroup|$(D_4)=G(8,3)$ means that $D_4$ is the third in
the list of groups of order $8$. We emphasize that all quasi-bundles
with $\chi(\mO_S)=1$ verify $K_S^2=8$ (see \cite{Se90}, Proposition
3.5), whereas imposing some RDPs allows us to obtain surfaces with
lower $K_S^2$. In particular, as a by-product of our classification,
we produce several examples with $p_g=q=1$ and $K_S^2=6$. In the
survey paper \cite{BaCaPi06} the minimal surfaces of general type
with these invariants are referred as ``mysterious''. Actually,
there was only one example hitherto known, described by C. Rito in
\cite{Ri07}. It verifies $g_{\textrm{alb}}=3$ and is obtained as a
double cover of a Kummer surface; the construction makes use of the
computer algebra program MAGMA in order to find a branch curve with
the right singularities. We note that Rito's surface is not a
standard isotrivial fibration, because the reducible fibres of its
Albanese pencil contain no $HJ$-strings (cf. Theorem \ref{Serrano}).
Therefore all examples with $p_g=q=1$, $K_S^2=6$ and
$g_{\textrm{alb}}=3$ described in the present paper were previously
unknown; in addition, we provide the first examples with
$g_{\textrm{alb}}=4$. Our viewpoint also sheds some new light on
surfaces with $p_g=q=1$, $K_S^2=4$ and $g_{\textrm{alb}}=2,3$. An
example with $K_S^2=4$ and
 $g_{\textrm{alb}}=2$ was previously given
 by Catanese (\cite{Ca99}) as the minimal resolution
 of a bidouble cover of $\mathbb{P}^2$;
examples with $K_S^2=4$ and $g_{\textrm{alb}}=3$ were constructed by
Ishida (\cite{Is05}) as the minimal resolution of a double cover
of the $2-$fold symmetric product $E^{(2)}$ of an elliptic curve.
Both covers of Catanese and Ishida contain non-rational
singularities, whereas in all our examples $T$ has only nodes; it
follows that all surfaces with $K_S^2=4$ presented here are new.
Finally, we obtain two examples with $K_S^2=2$; they can be also
constructed as double covers of $E^{(2)}$ and in
both case we describe the six-nodal branch curve in detail
(Proposition \ref{descriptionK2=2}). These two examples belong to
the same
 irreducible component of the moduli space of surfaces of general
 type with $K_S^2=2, \; \chi(\mO_S)=1$, which is in fact irreducible (\cite{Ca});
 then it would be desirable to know whether any two surfaces in our
list, with the same $K_S^2$ and $g_{\textrm{alb}}$, are deformation
equivalent. We conjecture that the answer is negative, but this
question is at the present not solved. One could obtain some partial
information by computing in every case the index of the
paracanonical system, which is a topological invariant
(\cite{CaCi91}, Theorem 1.4; see also \cite{Pol08}, Theorem 6.3),
but we will not develop this point here. \\
We shall now explain in more detail the steps of our classification
procedure. The crucial fact is that, since $G$ acts on both $C$ and
$F$, the geometry of $S$ is encoded in the geometry of the two
$G-$covers $h \colon C \lr C/G$, $f \colon F \lr F/G$. This allows
us to ``detopologize" the problem by transforming it into an
equivalent problem about the existence of a pair $(\mathcal{V},
\mathcal{W})$ of generating vectors for $G$ of type $(0 \; | \; m_1,
\ldots, m_r)$ and $(1 \; | \; n_1, \ldots, n_s)$, respectively (see
Section \ref{topological} for the definitions). These vectors must
satisfy some additional properties in order to obtain a quotient
$T=(C \times F)/G$ with only RDPs and whose desingularization $S$
has the desired invariants
(Proposition \ref{strukture}). \\
In Section \ref{topological} we present some preliminaries and we
fix the algebraic set-up. In Proposition \ref{riemann ext},
 which is essentially a reformulation of Riemann's existence theorem, we show
that a smooth projective curve $Y$ of genus $\mathfrak{g}'$ admits
 a $G-$cover $X \lr Y$, branched in $r$ points $P_1, \ldots, P_r$
 with branching numbers $m_1, \ldots, m_r$, if and only if
 $G$ contains a generating vector $\mathcal{V}$ of type
 $(\mathfrak{g}' \; | \; m_1, \ldots ,m_r)$.
For every $h \in G$ we give a formula that computes the number
of fixed points of $h$ on $X$ in terms of $\mathcal{V}$ (Proposition \ref{fixed-points}). \\
In Section \ref{non-generation} we collect some non-generation
results for finite groups  which will be useful in the sequel of the
paper. They are obtained either by direct computation or
 by using the GAP4 database of small groups.
For every group we refer to the presentation given
in the corresponding table of Appendix $A$. The reader that finds
these results too dry or boring  might skip this section for the
moment and come
back to it when reading Section \ref{general case}.  \\
In Section \ref{standard isotrivial} we establish the main properties of standard
isotrivial fibrations (following \cite{Se96})
and we compute their invariants in the case where $T$ has only RDPs. \\
In Sections \ref{chi=1} and \ref{structure theorem} we show that if $S$ is a
standard isotrivial fibration of general type with $p_g=q=1$ and $T$
contains only RDPs, then $S$ is a minimal surface, $K_S^2$ is even and the
singularities of $T$ are exactly $8-K_S^2$ nodes. Furthermore we prove Proposition
\ref{strukture}, which plays a crucial role in this paper as it
provides the translation
of our classification problem ``from geometry to algebra''. \\
In Section \ref{abelian case} we show our
Main Theorem assuming that the group $G$ is abelian; the proof is extended to
the nonabelian case in Section \ref{general case}. \\
The tables of Appendix $A$ contain the automorphism groups
 acting on Riemann surfaces of genus $2, \; 3$ and
$4$ so that the quotient is isomorphic to $\mathbb{P}^1$.
In the last two cases we listed only the
 nonabelian groups.
Tables \ref{2-abeliani}, \ref{2-nonabelian} and
 \ref{3-nonabelian} are adapted from [Br90, pages 252, 254, 255], whereas
Table \ref{4-nonabelian} is adapted from \cite[Theorem 1]{Ki03} and \cite{Vin00}.
For every $G$ we give a presentation, the branching data
and the \verb|IdSmallGroup|$(G)$.\\
Finally, in Appendix $B$ we give an example of GAP4
script used during the preparation of this work. \\

$\mathbf{Notations \; and \; conventions}$. All varieties,
morphisms, etc. in this article are defined over the field
$\mathbb{C}$ of the complex numbers. By ``surface'' we mean a
projective, non-singular surface $S$, and for such a surface $K_S$
denotes the canonical class, $p_g(S)=h^0(S, \; K_S)$ is the
\emph{geometric genus}, $q(S)=h^1(S, \; K_S)$ is the
\emph{irregularity} and $\chi(\mathcal{O}_S)=1-q(S)+p_g(S)$ is the
\emph{Euler characteristic}. If $T$ is a normal surface, a
 desingularization $\lambda \colon S \lr T$ is said to be \emph{minimal}
 if $\lambda$ does not contract any $(-1)-$curve in $S$. Such a minimal
 desingularization always exists and it is determined uniquely
 by $T$ (\cite{BPV}, p.86); it is worth pointing out that $S$ is not necessarily
 a minimal surface (cf. Proposition \ref{minimal}). \\
Throughout the paper we use the
following notation for groups:
\begin{itemize}
\item $\mZ_n$: cyclic group of order $n$.
\item $D_{p,q,r}=\mathbb{Z}_p \ltimes \mathbb{Z}_q= \langle x,y \; |
\; x^p=y^q=1, xyx^{-1}=y^r \rangle$: split metacyclic group of order
$pq$, note that $r^p \equiv 1 \; \; \textrm{mod}\; q$. The group
$D_{2,n,-1}$ is the dihedral group of order $2n$, that will be
denoted by $D_n$.
\item $S_n, \;A_n$: symmetric, alternating group on $n$ symbols. We write the
composition of permutations from the right to the left; for
instance, $(13)(12)=(123)$.
\item $\textrm{GL}_n(\mathbb{F}_q), \; \textrm{SL}_n(\mathbb{F}_q), \;
\textrm{PSL}_n(\mathbb{F}_ q)$: general linear, special linear and projective
special linear group of $n \times n$ matrices over a field with $q$
elements.
\item Whenever we give a presentation of a semi-direct product $H \ltimes N$, the
first generators represent $H$ and the last generators represent
$N$. The action of $H$ on $N$ is specified by conjugation relations.
\item The order of a finite group $G$ is denoted by $|G|$. If $H$ is
a subgroup of $G$, the centralizer of $H$ in $G$ is denoted by
$C_G(H)$ and the normalizer of $H$ in $G$ by $N_G(H)$. The conjugacy
relation in $G$ is denoted by $\sim_G$.
\item The subgroup generated by $x_1, \ldots, x_n \in G$
 is denoted by $\langle x_1, \ldots, x_n \rangle$. The derived subgroup of $G$ is denoted by $G'$.
 The center of $G$ is denoted by $Z(G)$.
 The set of elements of $G$ different from the identity is denoted
 by $G^{\times}$.
\item If $x \in G$, the order of $x$ is denoted by $o(x)$ and
the conjugacy class of $x$ by $\textrm{Cl}(x)$. If $x,y \in G$,
their commutator is defined as $[x,y]=xyx^{-1}y^{-1}$.
\item All groups are represented in multiplicative format. \\
\end{itemize}
$\mathbf{Acknowledgements.}$ The author wishes to thank I. Bauer, A.
Broughton, G. Carnovale, F. Catanese, C. Ciliberto, G. Infante, J.
Van Bon for useful discussions and suggestions. Part of this
research was done during the winter semester of the academic year
2006-2007, when the
 author was visiting the Department of Mathematics at the Imperial College
(London, UK) and was supported by a ``Accademia dei Lincei" grant.
He is especially grateful to A. Corti for his warm hospitality.

\section{Algebraic background} \label{topological}
In this section we present some preliminaries and we fix the
algebraic set-up. Many of the result that we collect here are
standard, so proofs are often omitted. We refer the reader to [Br90,
Section 2], [Bre00, Chapter 3], \cite{H71} and [Pol08, Section 1]
for more details.
\begin{definition} \label{generating vect}
Let $G$ be a finite group and let
\begin{equation*}
\mathfrak{g}' \geq 0, \quad   m_r \geq m_{r-1} \geq \ldots \geq m_1
\geq 2
\end{equation*}
be integers. A \emph{generating vector} for $G$ of type
$(\mathfrak{g}' \; | \; m_1, \ldots ,m_r)$ is a $(2
\mathfrak{g}'+r)$-tuple of elements
\begin{equation*}
\mathcal{V}=\{g_1, \ldots, g_r; \; h_1, \ldots, h_{2\mathfrak{g}'}
\}
\end{equation*}
such that the following conditions are satisfied:
\begin{itemize}
\item the set $\mathcal{V}$ generates $G$;
\item $o(g_i)=m_i$;
\item $g_1g_2\cdots g_r \Pi_{i=1}^{\mathfrak{g}'} [h_i,h_{i+\mathfrak{g}'}]=1$.
\end{itemize}
If such a $\mathcal{V}$ exists, then $G$ is said to be
$(\mathfrak{g}' \; | \; m_1, \ldots ,m_r)-$\emph{generated}.
\end{definition}
For convenience we make abbreviations such as $(4 \;| \; 2^3, 3^2)$
for $(4 \; | \; 2,2,2,3,3)$ when we write down the type of the
generating vector $\mathcal{V}$.
\begin{proposition} \label{ab-no1}
If an abelian group $G$ is $(\mathfrak{g}'\; |\; m_1, \ldots, m_r)-$generated,
then $r \neq 1$.
\end{proposition}
\begin{proof}
If $r=1$ and $\mathcal{V}=\{g_1, \; h_1, \ldots, h_{2\mathfrak{g}'}
\}$ is a generating vector, we have
\begin{equation*}
1=g_1\Pi_{i=1}^{\mathfrak{g}'} [h_i,h_{i+\mathfrak{g}'}]=g_1,
\end{equation*}
a contradiction because $o(g_1)=m_1 \geq 2$.
\end{proof}
The following result, which is essentially a reformulation of
Riemann's existence theorem, translates the problem of finding Riemann
surfaces with automorphisms into the group theoretic problem of
finding groups $G$ which contain suitable generating vectors.
\begin{proposition} \label{riemann ext}
A finite group $G$ acts as a group of automorphisms of some compact
Riemann surface $X$ of genus $\mathfrak{g}$ if and only if there
exist integers $\mathfrak{g}' \geq 0$ and  $m_r \geq m_{r-1} \geq
\ldots \geq m_1 \geq 2$ such that $G$ is $(\mathfrak{g}'\; |\; m_1,
\ldots, m_r)-$generated, with generating vector $\mathcal{V}=\{g_1,
\ldots, g_r; \; h_1, \ldots, h_{2\mathfrak{g}'} \}$, and the
following Riemann-Hurwitz relation holds:
\begin{equation} \label{riemanhur}
2\mathfrak{g}-2=|G| \left(
2\mathfrak{g}'-2+\sum_{i=1}^r\bigg(1-\frac{1}{\;m_i} \bigg)
 \right).
\end{equation}
If this is the case then $\mathfrak{g}'$ is the genus of the
quotient Riemann surface $Y:=X/G$ and the $G-$cover $X \lr Y$ is
branched in $r$ points $P_1, \ldots, P_r$ with branching numbers
$m_1, \ldots, m_r$, respectively. In addition, the subgroups
$\langle g_i \rangle$ and their conjugates provide all the nontrivial
stabilizers of the action of $G$ on $X$.
\end{proposition}
Let $G$, $\mathcal{V}$ and $X$ be as in Proposition \ref{riemann ext}.
For any $h \in G$ set $H:=\langle h \rangle$ and define
\begin{equation*}
\textrm{Fix}_X(h)=\textrm{Fix}_X(H):=\{x \in X \; |\; hx=x \}.
\end{equation*}

\begin{proposition} \label{fixed-points}
If $o(h)=m$ then
\begin{equation*} %\label{formula-per-fix}
|\emph{Fix}_X(h)|=|N_G(H)| \cdot \sum_{\substack{1 \leq i \leq r \\
m|m_i
\\ H \sim_G \; \langle g_i^{m_i/m} \rangle }} \frac{1}{\; m_i}.
\end{equation*}
\end{proposition}
\begin{proof}
(See \cite{Bre00}, Lemma 10.4). Let $x$ be in $\textrm{Fix}_X(h)$ and let $R_i$
be a set of coset representatives of $\langle g_i \rangle$ in $G$. Then
\begin{equation*}
\begin{split}
\textrm{Fix}_X(h)&=\biguplus_{1 \leq i \leq r} \{ \sigma x \; | \; \sigma \in R_i,
\; H \leq \langle \sigma g_i \sigma^{-1} \rangle \} \\
        &=\biguplus_{1 \leq i \leq r} \{ \sigma x \; | \; \sigma \in R_i,
\; H = \langle \sigma g_i^{m_i/m} \sigma^{-1} \rangle \}.
\end{split}
\end{equation*}
Taking the cardinalities on both sides, we get
\begin{equation*}
\begin{split}
|\textrm{Fix}_X(h)|&=\sum_{1 \leq i \leq r} | \{ \sigma x \; | \; \sigma \in R_i, \; H = \langle \sigma g_i^{m_i/m} \sigma^{-1} \rangle \} | \\
&=\sum_{1 \leq i \leq r} | \{ \sigma \in R_i \;|\;  H = \langle \sigma g_i^{m_i/m} \sigma^{-1} \rangle \} | \\
&=\sum_{1 \leq i \leq r} \frac{1}{m_i} | \{ \sigma \in G \; | \; H = \langle \sigma g_i^{m_i/m} \sigma^{-1} \rangle \} |,
\end{split}
\end{equation*}
where the set in the $i-$th summand has cardinality $|N_G(H)|$ if $H$ is
$G-$conjugate to $\langle g_i^{m_i/m} \rangle$, and is empty otherwise.
\end{proof}

\begin{corollary} \label{ordine2}
If $o(h)=2$ then
\begin{equation} \label{formula-per-fix-ord2}
|\emph{Fix}_X(h)|=\frac{|G|}{|\emph{Cl}(h)|} \cdot
\sum_{\substack{2|m_i
\\ H \sim_G \; \langle g_i^{m_i/2} \rangle }} \frac{1}{\; m_i}.
\end{equation}
If $o(h)=2$ and $h \in Z(G)$ then
\begin{equation} \label{formula-per-fix-ab}
|\emph{Fix}_X(h)|=|G| \cdot \sum_{\{i \, | \, h \in \langle g_i
\rangle \}} \frac{1}{\; m_i}.
\end{equation}
\end{corollary}
\begin{proof}
Since $H \cong \mathbb{Z}_2$ we have $N_G(H)=C_G(H)$, so Proposition
\ref{fixed-points} implies (\ref{formula-per-fix-ord2}). The proof
of (\ref{formula-per-fix-ab}) is now immediate.
\end{proof}

\section{Some non-generation results} \label{non-generation}

This section contains some non-generation results for finite groups
which will be useful in the sequel of our classification procedure.
They are obtained either by direct computation or by using the GAP4
database of small groups. We will first use them in Section
\ref{general case}. For every group we refer to the presentation
given in the corresponding table of Appendix $A$.

\begin{lemma} \label{no 1-22}
Let $G$ be a nonabelian finite group containing a unique element
$\ell$ of order $2$. Then $G$ is not $(1\; | \; 2^2)-$generated.
\end{lemma}
\begin{proof}
Assume that $G$ is $(1\; | \; 2^2)-$generated, with generating
vector $\mathcal{V}= \{\ell_1, \ell_2; \; h_1,h_2\}$. Since $\ell$
is the only element of order $2$ in $G$, it follows $\ell \in Z(G)$ and
 $\ell_1=\ell_2=\ell$, hence $[h_1, h_2]=1$. Therefore $G=\langle
\ell, h_1, h_2 \rangle$ would be abelian, a contradiction.
\end{proof}

\begin{proposition} \label{no-generation-1-22}
Referring to Table \emph{\ref{2-nonabelian}} of Appendix A, the groups $G$
in cases $(2b)$, $(2d)$, $(2h)$ are not $(1\; | \; 2^2)-$generated.
\end{proposition}
\begin{proof}
It is sufficient to show that $G$ satisfies the hypotheses of Lemma \ref{no 1-22}. \\
$\bullet$ Case ($2b$). $G=Q_8$. Take $\ell=-1$.
\\
$\bullet$ Case ($2d$). $G=D_{4,3,-1}$. Take $\ell=x^2$. \\
$\bullet$ Case($2h$). $G=\textrm{SL}_2(\mathbb{F}_3)$.
Take $\ell=\left( \begin{array}{cc}
    -1 & \;\;0 \\
    \;\;0 & -1 \\
  \end{array}
\right)$.
\end{proof}

\begin{proposition} \label{no-generation-1-21}
Referring to Table \emph{\ref{2-nonabelian}} of Appendix A, the groups $G$ in cases
$(2d)$, $(2e)$, $(2f)$, $(2g)$,
$(2h)$, $(2i)$  are not  $(1\; | \; 2^1)-$generated.
\end{proposition}
\begin{proof}
We do a case-by-case analysis. \\ \\
$\bullet$ Case $(2d)$. $G=D_{4,3,-1}$.\\
Looking at the presentation of $G$, one checks that $G'=\langle y \rangle \cong \mathbb{Z}_3$. Therefore $G$ contains no commutators of
order $2$, so it cannot be $(1\; | \; 2^1)-$generated.  \\
\\
$\bullet$ Case $(2e)$. $G=D_6$.\\
We have $G'=\langle y^2 \rangle \cong \mathbb{Z}_3$, so $G$ contains
no commutators of order $2$. \\ \\
$\bullet$ Case $(2f)$. $G=D_{2,8,3}$. \\
We have $G'= \langle y^2 \rangle \cong \mathbb{Z}_4$ and the only
commutator of order $2$ is $y^4$. A direct computation shows that if
$[h_1,h_2]=y^4$ then either $\langle h_1, h_2 \rangle \cong D_4$ or
$\langle h_1, h_2 \rangle \cong Q_8$. In
particular $\langle h_1, h_2 \rangle \neq G$, hence $G$ is
not  $(1\; | \; 2^1)-$generated.\\ \\
$\bullet$ Case $(2g)$. $G=\mathbb{Z}_2 \ltimes ((\mathbb{Z}_2)^2
\times \mathbb{Z}_3)=G(24,8)$. \\
We have $G'= \langle y,w \rangle \cong \mathbb{Z}_6$ and the only commutator of order $2$ is $y$. If $[h_1, h_2]=y$ then $\langle
h_1, h_2 \rangle \cong
D_4$, so $G$ is not  $(1\; | \; 2^1)-$generated.\\ \\
$\bullet$ Case $(2h)$. $G=\textrm{SL}_2(\mathbb{F}_3)$. \\
The group $G$ contains a unique element
of order $2$, namely
$\ell=     \left(
              \begin{array}{cc}
                -1 & \;\;0 \\
                \;\;0 & -1 \\
              \end{array}
            \right)$.
A direct computation shows that $G'\cong Q_8$ and that $\ell$
can be expressed as a commutator in $24$ different ways. Moreover,
if $[h_1,h_2]=\ell$ we have  $\langle h_1,h_2 \rangle \cong Q_8$,
so $G$ is not  $(1\; | \; 2^1)-$generated.\\ \\
$\bullet$ Case $(2i)$. $G=\textrm{GL}_2(\mathbb{F}_3)$. \\
It is well known that $G'=\textrm{SL}_2(\mathbb{F}_3)$; then $G'$
contains a unique element of order $2$, namely $\ell$. Either by
direct computation or by using GAP4, one can check that there are
$96$ different ways to write $\ell$ as a commutator in $G$. If
$[h_1,h_2]=\ell$ and both  $h_1$ and $h_2$ belong to
$\textrm{SL}_2(\mathbb{F}_3)$, then $\langle h_1, h_2 \rangle \cong
Q_8$; otherwise $\langle h_1, h_2 \rangle \cong D_4$. In both cases
$\langle h_1, h_2 \rangle \neq G$, hence $G$ is not  $(1\; | \;
2^1)-$generated.
\end{proof}

\begin{proposition} \label{Table-3-no-1-2-gen}
Referring to Table \emph{\ref{3-nonabelian}}  of Appendix A, the groups $G$ in cases $(3d)$, $(3e)$,
$(3i)$, $(3j)$, $(3l)$, $(3n)$, $(3o)$, $(3p)$, $(3q)$,
$(3r)$, $(3s)$, $(3t)$, $(3u)$, $(3v)$, $(3w)$ are not $(1\; | \;2^1)-$generated.
\end{proposition}
\begin{proof}
We have already proven the statement in cases $(3d)$, $(3e)$ and
$(3n)$: see Proposition
\ref{no-generation-1-21}, cases  $(2d)$, $(2e)$ and $(2h)$.
Now let us consider the remaining
cases. \\ \\
$\bullet$ Case $(3i)$. $G=\mathbb{Z}_2 \times D_4$. \\
The group $G$ cannot be generated by two elements, so in particular
it cannot be $(1\; | \;2^1)-$generated.
 \\ \\
$\bullet$ Case $(3j)$. $G=\mathbb{Z}_2 \ltimes(\mathbb{Z}_2 \times
\mathbb{Z}_4)=G(16,13)$.\\
We have $G'= \langle z^2 \rangle \cong \mathbb{Z}_2$. By direct
computation or by using GAP4 (see Appendix $B$ for the corresponding script)
 we can check that if $[h_1, h_2]=z^2$
then either $\langle h_1, h_2 \rangle \cong D_4$ or
$\langle h_1, h_2 \rangle \cong Q_8$, so $\langle h_1, h_2
\rangle \neq G$. \\ \\
$\bullet$ Case $(3l)$. $G=D_{2,12,5}$. \\
We have $G'=\langle y^4 \rangle \cong \mathbb{Z}_3$, so $G$ contains
no commutators of order $2$. \\ \\
$\bullet$ Cases $(3o)$ and $(3p)$. $G=S_4$. \\
We have $G'=A_4$. If $[h_1, h_2]$ has order $2$ then $\langle h_1,
h_2 \rangle \cong D_4$  or $\langle h_1, h_2 \rangle \cong A_4$, so
 $G$ is not $(1\; | \;2^1)-$generated. \\ \\
$\bullet$ Case $(3q)$. $G=\mathbb{Z}_2 \ltimes (\mathbb{Z}_2 \times
\mathbb{Z}_8)=G(32,9)$. \\
We have $G'=\langle yz^2 \rangle \cong \mathbb{Z}_4$ and the only
commutator of order $2$ is $(yz^2)^2=z^4$. If $[h_1, h_2]=z^4$ then
$\langle h_1, h_2 \rangle$ has order $8$ or $16$, hence $\langle
h_1, h_2 \rangle \neq G$. \\ \\
$\bullet$ Case $(3r)$. $G=\mathbb{Z}_2 \ltimes D_{2,8,5}=G(32,11)$. \\
We have $G'=\langle yz^2 \rangle \cong \mathbb{Z}_4$ and the only
commutator of order $2$ is $(yz^2)^2=z^4$. If $[h_1, h_2]=z^4$ then
$\langle h_1, h_2 \rangle$ has order $8$ or $16$, hence $\langle
h_1, h_2
\rangle \neq G$. \\ \\
$\bullet$ Case $(3s)$. $G=\mathbb{Z}_2 \times S_4$. \\
If $[h_1, h_2]$ has order $2$ then $|\langle h_1, h_2
\rangle| \leq 24$, so $\langle h_1, h_2
\rangle \neq G$.    \\ \\
$\bullet$ Case $(3t)$. $G=G(48,33)$. \\
We have $G' = \langle t,z, w \rangle \cong Q_8$ and the only
commutator of order $2$ is $t$. If $[h_1, h_2]=t$ then $\langle h_1,
h_2 \rangle
 \cong D_4$ or $\langle h_1, h_2 \rangle \cong Q_8$, so $G$ is not $(1\; | \;2^1)-$generated. \\ \\
$\bullet$ Case $(3u)$. $G=\mathbb{Z}_3 \ltimes (\mathbb{Z}_4)^2=G(48,3)$.\\
We have $G'=\langle y, z \rangle \cong (\mathbb{Z}_4)^2$. If $[h_1,
h_2]$ has order $2$ then $\langle h_1, h_2
\rangle \cong A_4$, so $G$ is not $(1\; | \;2^1)-$generated. \\ \\
$\bullet$ Case $(3v)$. $G=S_3 \ltimes (\mathbb{Z}_4)^2=G(96,64)$. \\
We have $G'= \langle y,z \rangle$ and  $|G'|=48$. The elements of
order $2$ in $G'$ are $z^2$, $y^2z^2y$, $yz^2y^2$. If $[h_1, h_2]$
has order $2$ then $|\langle h_1, h_2 \rangle| \leq 16$, so $G$ is not $(1\; | \;2^1)-$generated. \\ \\
$\bullet$ Case $(3w)$. $G=\textrm{PSL}_2(\mathbb{F}_7)$. \\
Since $G$ is simple we have $G'=G$. If $[h_1, h_2]$ has order $2$
then either $\langle h_1, h_2 \rangle
\cong D_4$ or $\langle h_1, h_2 \rangle \cong A_4$, so  $G$ is not $(1\; | \;2^1)-$generated. \\ \\
\end{proof}
\begin{proposition} \label{Table-3-no-1-4-gen}
Referring to Table \emph{\ref{3-nonabelian}} of Appendix A, the
groups $G$ in cases $(3i)$,
$(3j)$, $(3s)$, $(3v)$ are not $(1\; |
\;4^1)-$generated.
\end{proposition}
\begin{proof}
We do a case-by-case analysis.\\  \\
$\bullet$ Case $(3i)$. $G=\mathbb{Z}_2 \times D_4$. \\
We have $G'=\langle (1, y^2) \rangle \cong \mathbb{Z}_2$; therefore
$G$ contains no commutators of order $4$ and so it cannot be $(1\; |
\;4^1)-$generated.\\ \\
$\bullet$ Case $(3j)$. $G=\mathbb{Z}_2 \ltimes (\mathbb{Z}_2 \times
\mathbb{Z}_4)=G(16,13)$.\\
We have $G'=\langle z^2 \rangle \cong \mathbb{Z}_2$, so $G$ contains
no commutators of order $4$ and we conclude as in the previous case.  \\ \\
$\bullet$ Case $(3s)$. $G=\mathbb{Z}_2 \times S_4$. \\
We have $G' \cong A_4$, so $G$ contains no commutators of order
$4$. \\ \\
$\bullet$ Case $(3v)$. $G=S_3 \ltimes (\mathbb{Z}_4)^2=G(96,64)$. \\
If $[h_1, h_2]$ has order $4$ then $|\langle h_1, h_2
\rangle| \leq 48$, so $G$ is not  $(1\; |
\;4^1)-$generated.
\end{proof}

\begin{proposition} \label{Table-4-no-1-2-gen}
Referring to  Table \emph{\ref{4-nonabelian}} of Appendix A, the
groups $G$ in cases $(4g)$, $(4h)$, $(4i)$,
 $(4j)$, $(4k)$, $(4l)$, $(4o)$, $(4p)$, $(4q)$, $(4s)$, $(4t)$, $(4u)$,
 $(4v)$, $(4w)$, $(4y)$, $(4z)$, $(4aa)$, $(4ab)$ are not $(1\; | \;2^1)-$generated.
\end{proposition}
\begin{proof}
Again a case-by-case analysis. \\ \\
$\bullet$ Cases $(4g)$ and $(4h)$. $G=D_6$. \\
See Proposition \ref{no-generation-1-21}, case $(2e)$. \\ \\
$\bullet$ Case $(4i)$. $G=D_8$. \\
We have $G'= \langle y^2 \rangle \cong \mZ_4$ and the only commutator
of order $2$ is $y^4$. If $[h_1, h_2]=y^4$ then $\langle h_1, h_2
\rangle \cong D_4$, hence
$G$ is not $(1\; | \;2^1)-$generated. \\ \\
$\bullet$ Case $(4j)$. $G=G(16,9)$. \\
We have $G'= \langle z \rangle \cong \mZ_4$ and the only commutator
of order $2$ is $z^2$. If $[h_1, h_2]=z^2$ then $\langle h_1, h_2
\rangle \cong Q_8$, hence
$G$ is not $(1\; | \;2^1)-$generated. \\ \\
$\bullet$ Cases $(4k)$ and $(4l)$. $G=\mZ_3 \times S_3$. \\
We have $G' \cong \mZ_3$, so $G$ contains no commutators of order $2$. \\
\\
$\bullet$ Case $(4o)$. $G=D_{4,5,-1}$. \\
We have $G'= \langle y \rangle \cong \mZ_5$, so $G$ contains no commutators of order $2$. \\ \\
$\bullet$ Case $(4p)$. $G=D_{4,5,2}$. \\
We have $G'= \langle y \rangle \cong \mZ_5$, so $G$ contains no commutators of order $2$. \\ \\
$\bullet$ Case $(4q)$. $G=S_4$. \\
See Proposition \ref{Table-3-no-1-2-gen}, cases $(3o)$ and $(3p)$.  \\ \\
$\bullet$ Case $(4s)$. $G=\textrm{SL}_2(\mathbb{F}_3)$. \\
See Proposition \ref{no-generation-1-21}, case $(2h)$. \\ \\
$\bullet$ Case $(4t)$. $G=D_{2,16,7}$. \\
We have $G'= \langle y^2 \rangle \cong \mZ_8$ and the only commutator of
order $2$ is $y^8$. If $[h_1, h_2]=y^8$ then $\langle h_1, h_2 \rangle \cong D_4$ or $\langle h_1, h_2 \rangle \cong Q_8$, hence $G$ is not $(1\; | \;2^1)-$generated. \\ \\
$\bullet$ Cases $(4u)$ and $(4v)$. $G=(\mZ_2)^2 \ltimes (\mZ_3)^2=G(36,10)$ \\
We have $G'= \langle z,w \rangle \cong \mZ_3 \times \mZ_3$, so $G$ contains
no commutators of order $2$. \\ \\
$\bullet$ Case $(4w)$. $G=\mZ_6 \times S_3$. \\
We have $G' \cong \mZ_3$, so $G$ contains
no commutators of order $2$. \\ \\
$\bullet$ Case $(4y)$. $G=\mZ_4 \ltimes (\mZ_3)^2=G(36,9)$. \\
We have $G'= \langle y,z \rangle \cong \mZ_3 \times \mZ_3$, so $G$ contains
no commutators of order $2$. \\ \\
$\bullet$ Case $(4z)$. $G=D_4 \ltimes \mZ_5=G(40,8)$. \\
We have $G'= \langle y^2,z \rangle \cong \mZ_{10}$ and the only commutator
of order $2$ is $y^2$. If $[h_1, h_2]=y^2$ then $\langle h_1, h_2 \rangle \cong D_4$, hence $G$ is not $(1\; | \;2^1)-$generated. \\ \\
$\bullet$ Case $(4aa)$. $G=A_5$. \\
Since $G$ is simple we have $G'=G$. If $[h_1, h_2]$ has order $2$ then $\langle h_1, h_2 \rangle \cong A_4$, hence $G$ is not $(1\; | \;2^1)-$generated. \\ \\
$\bullet$ Case $(4ab)$. $G=\mZ_3 \times S_4$ \\
We have $G'= A_4$. If $[h_1, h_2]$ has order $2$ then $|\langle h_1, h_2 \rangle| \leq 36$, hence $G$ is not $(1\; | \;2^1)-$generated. \\ \\
\end{proof}

\section{Standard isotrivial fibrations} \label{standard isotrivial}
In this section we establish the basic properties of standard
isotrivial fibrations. Definition \ref{def-stand} and Theorem \ref{Serrano}
 can be found in \cite{Se96}.
\begin{center}
\emph{From now on, $S$ will always denote a smooth, projective surface of general type}.
\end{center}
\begin{definition} \label{def-stand}
We say that $S$ is a \emph{standard isotrivial fibration} if
there exists a finite group $G$ acting faithfully on two smooth
projective curves $C$ and $F$ so that $S$ is isomorphic to the
minimal desingularization of $T:=(C \times F)/G$. The two maps
$\alpha \colon S \lr C/G$, $\beta \colon S \lr F/G$ will be referred
as the \emph{natural projections}. If $T$ is smooth then $S=T$ is
called a $\emph{quasi-bundle}$, or a
 surface $\emph{isogenous to an unmixed product}$.
\end{definition}

The stabilizer $H  \subseteq G$ of a point $y \in F$ is a cyclic
group (\cite{FK92}, p.106). If $H$ acts freely on $C$, then $T$ is
smooth along the scheme-theoretic fibre of $\sigma \colon T \lr F/G$
 over $\bar{y} \in F/G$, and this fibre consists of the curve $C/H$
 counted with multiplicity $|H|$. Thus, the smooth fibres of $\sigma$
 are all isomorphic to $C$. On the contrary, if $ x \in C$ is fixed
 by some non-zero element of $H$, then $T$ has a cyclic quotient
 singularity  over the point $\overline{(x,y)} \in (C \times F)/G$. In
 this case, the fibre of $\overline{(x,y)}$ on the minimal
 desingularization $\lambda \colon S \lr T$ is an $HJ$-string
 (abbreviation of Hirzebruch-Jung string), that is to say, a
 connected union of smooth rational curves $Z_1, \ldots, Z_n$ with
 self-intersection $\leq 2$, and ordered linearly so that $Z_i
 Z_{i+1}=1$ for all $i$, and $Z_iZ_j=0$ if $|i-j| \geq 2$
 (\cite{BPV}, III 5.4). These observations lead to the following
 statement, which describes the singular fibres that can arise in a
 standard isotrivial fibration (see \cite{Se96}, Theorem 2.1).
\begin{theorem} \label{Serrano}
Let $\lambda \colon S \lr T=(C \times F)/G$ be a standard
isotrivial fibration and let us consider the natural projection
 $\beta \colon S \lr F/G$.
Take any point over $\bar{y} \in F/G$ and let $\Lambda$ denote the
fibre of $\beta$ over $\bar{y}$. Then
\begin{itemize}
\item[$(i)$] the reduced structure of $\Lambda$ is the union of an
irreducible curve $Y$, called the central component of $\Lambda$,
 and either none or at least two mutually disjoint HJ-strings, each
 meeting $Y$ at one point. These strings are in one-to-one
 correspondence with the branch points of $C \lr C/H$, where $H
 \subseteq G$ is the stabilizer of $y$;
\item[$(ii)$] the intersection of a string with $Y$ is transversal,
and it takes place at only one of the end components of the string;
\item[$(iii)$] $Y$ is isomorphic to $C/H$, and has multiplicity
equal to $|H|$ in $\Lambda$.
\end{itemize}
Evidently, a completely similar statement holds if we consider the
natural projection $\alpha \colon S \lr C/G$.
\end{theorem}

\begin{remark} \label{An-cycles}
The $HJ$-strings arising from the minimal resolution of RDPs are
precisely the $A_n-$cycles.
\end{remark}

Theorem \ref{Serrano} and Remark \ref{An-cycles} now imply

\begin{corollary} \label{singularities-0}
Let us suppose that $T$ has at worst \emph{RDPs}, and let $\Lambda$
be any fiber of $\beta \colon S \lr F/G$. Then $\Lambda$ contains
either none or at least two $A_n-$cycles. An analogous statements
holds if we consider any fibre $\Phi$ of $\alpha \colon S \lr C/G$.
\end{corollary}

It is worth pointing out that a standard isotrivial fibration is not
necessarily a minimal surface; indeed, the central component of some
reducible fibre might be a $(-1)-$curve. A criterion for minimality
is provided by the following

\begin{proposition} \label{minimal}
If $T$ has at worst \emph{RDPs} then both fibrations $\alpha \colon S \lr C/G$
and $\beta \colon S \lr F/G$ are relatively minimal. In addition,
if either $g(C/G)>0$ or $g(F/G)>0$ then $S$ is a minimal model.
\end{proposition}
\begin{proof}
Let us suppose that  $\beta$ is not relatively minimal;
then there is a singular fibre $\Lambda$ whose
central component $Y$ is a $(-1)-$curve. Corollary \ref{singularities-0} implies that $\Lambda$
 contains (at least) two disjoint $A_n-$cycles $Z_1$, $Z_2$ such that $YZ_1=YZ_2=1$. Thus
 by blowing down $Y$ we obtain a surface $S'$ with two $(-1)-$curves $E_1$,
 $E_2$ such that $E_1E_2=1$, a contradiction because
 $S$ is of general type (cf. \cite{BPV}, Proposition 4.6 p.79).
The proof for $\alpha$ is similar.
The last part of the statement follows at once because a fibration over a curve
 of strictly positive genus is minimal if and only if it is relatively minimal.
\end{proof}

Now set $\mathfrak{g}_1':=g(F/G)$ and $\mathfrak{g}_2':=g(C/G)$. By
Proposition \ref{riemann ext} it follows that there exist
\begin{itemize}
\item[-] integers $2 \leq m_1 \leq m_2 \leq \ldots \leq m_r$ such that $G$ is
$(\mathfrak{g}_1' \; | \; m_1 , \ldots, m_r)-$generated and
\item[-] integers $2 \leq n_1 \leq n_2 \leq \ldots \leq n_s$ such that $G$ is
$(\mathfrak{g}_2' \; | \; n_1 , \ldots, n_s)-$generated.
\end{itemize}

\begin{proposition} \label{ser:1}
If $T$ has at worst \emph{RDPs}, then
\begin{itemize}
\item $m_i\;$ divides $\;2g(C)-2$ \quad for all $\;i \in \{1,
  \ldots, r \}$;
\item $n_j\;$ divides $\;2g(F)-2$ \quad for all $\;j \in
  \{1, \ldots, s \}$.
\end{itemize}
\end{proposition}
\begin{proof}
Take any $\;i \in \{1,
  \ldots, r \}$. By Theorem \ref{Serrano} there exists a fibre
$\Lambda$ of $\beta \colon S \lr F/G$ having the form
$\Lambda=Y+Z$, where $Y$ is a component of multiplicity $m_i$ and
$Z$ is a (possibly empty) union of $(-2)-$curves. Setting $Y=m_iY'$
we obtain $K_S \Lambda =m_iK_SY'$; since $\Lambda$ is
algebraically equivalent to $C$ this implies $2g(C)-2=m_iK_SY'$. Thus
 $m_i$ divides $2g(C)-2$. Clearly, we can prove the second
  claim in the same way.
\end{proof}
\begin{corollary} \label{ser:2}
Assuming that $T$ has at worst \emph{RDPs}, the following holds:
\begin{itemize}
\item if $g(C)=2\;$ then $\;m_i=2$ \quad  for all $\;i \in \{1,
  \ldots, r \}$;
\item if $g(F)=2\;$ then $\;n_j=2$ \quad  for all $\;j \in \{1,
  \ldots, s \}$.
\end{itemize}
\end{corollary}
\begin{corollary} \label{nodal}
Suppose that $T$ has only \emph{RDPs}. If either $g(C)=2$ or $g(F)=2$ then
$T$ has at worst nodes $($i.e. $t_n=0$ for $n \geq 2).$
\end{corollary}
\begin{proof}
If $g(C)=2$ then by Corollary \ref{ser:2} it follows that the
non-trivial stabilizers of the action of $G$ on $F$ are isomorphic
to $\mathbb{Z}_2$, and this implies that the singularities of $T$
are at worst nodes. If $g(F)=2$  the argument is the same.
\end{proof}

The invariants of $S$ can be computed using
\begin{proposition} \label{invariants}
Let $V$ be a smooth algebraic surface, and let $G$ be a finite group
acting on $V$ with only isolated fixed points. Suppose that the
quotient $T:=V/G$ has at worst \emph{RDPs}, and let $\lambda \colon
S \lr T$ be the minimal desingularization. Let $t_n$ be the number
of singular points of type $A_n$ in $T$. Then we have
\begin{itemize}
\item[$(i)$] $|G| \cdot K_S^2=K_V^2$.
\item[$(ii)$] $|G| \cdot e(S)= e(V)+|G| \cdot \sum_n \frac{(n+1)^2-1}{n+1}t_n$.
\item[$(iii)$] $H^0(S, \Omega_S^1)=H^0(V, \Omega_V^1)^G$.
\end{itemize}
\end{proposition}
\begin{proof}
$(i)$ This is immediate because $G$ acts on $V$ with only isolated
fixed points and the singularities of $T$ are at worst RDPs.\\
$(ii)$ Let $\pi \colon V \lr T$ be the projection, $T^o$ be the
smooth locus of $T$ and $V^o:=\pi^{-1}(T^o)$; finally set
$S^o=\lambda^{-1}(T^o)$. Let $p \in T$ be a singularity of type
$A_n$; since $p$ is covered by $\frac{|G|}{n+1}$ points in $V$,
 we obtain
\begin{equation*}
e(V^o)=e(V)-\sum_n \frac{|G|}{n+1}t_n.
\end{equation*}
On the other hand, since $G$ acts on $V^o$ without fixed points, we
have $|G| \cdot e(S^o)=e(V^o)$. Finally, notice that $S$ is
 obtained from $S^o$ by attaching all the $A_n$-cycles; since every
 $A_n$-cycle $Z$ verifies $e(Z)=n+1$, the additivity of the Euler
 number implies
\begin{equation*}
\begin{split}
|G| \cdot e(S)&=|G| \cdot e(S^o)+|G| \cdot \sum_n (n+1)t_n \\
               &= e(V)+|G| \cdot \sum_n \frac{(n+1)^2-1}{n+1}t_n.
\end{split}
\end{equation*}
$(iii)$ See \cite{Fre71}.
\end{proof}
So we have
\begin{proposition} \label{invariantsS}
Let $\lambda \colon S \lr T=(C \times F)/G$ be a standard isotrivial
fibration such that $T$ has at worst \emph{RDPs}. Denote by
$t_n$ the number of singular points of type $A_n$ in $T$. Then the
invariants of $S$ are
\begin{itemize}
\item $K_S^2= \frac{8(g(C)-1)(g(F)-1)}{|G|}$
\item $e(S)=\frac{4(g(C)-1)(g(F)-1)}{|G|} + \sum_n
\frac{(n+1)^2-1}{n+1}t_n$
\item $q(S)=g(C/G)+g(F/G)$.
\end{itemize}
\end{proposition}
In particular this implies (cf. \cite{Se96}):
\begin{corollary}
The following are equivalent:
\begin{itemize}
\item $t_n=0$ for any $n \geq 1$ ;
\item $K_S^2=2e(S)$;
\item $S$ is a quasi-bundle.
\end{itemize}
\end{corollary}

\begin{remark} \label{singularities-2}
By Corollary \ref{singularities-0} it follows $\sum_{n} t_n \neq 1$.
\end{remark}

\section{The case $\chi(\mO_S)=1$} \label{chi=1}

\begin{proposition} \label{poss-sing}
Let $\lambda \colon S \lr T=(C \times F)/G$ be a standard isotrivial fibration,
such that $T$ contains at worst \emph{RDPs}.
In addition, let us assume $\chi(\mO_S)=1$. Then there are the following possibilities:
\begin{itemize}
\item $1 \leq K_S^2 \leq 8$ and $T$ contains $8-K_S^2$ points of type
$A_1$;
\item $K_S^2=3$ and $T$ contains two points of type $A_3$;
\item $K_S^2=2$ and T contains one point of type $A_1$ and two
points of type $A_3$;
\item $K_S^2=1$ and $T$ contains two points of type $A_1$ and two
points of type $A_3$.
\end{itemize}
\end{proposition}
\begin{proof}
If a minimal surface of general type with $\chi(\mO_S)=1$ contains
 some $A_n-$cycle then $n \leq 10$ (\cite{Mi}). Thus by
Proposition \ref{invariantsS} we have
\begin{equation*}
\begin{split}
&\frac{1}{2}K_S^2+\frac{3}{2}t_1+\frac{8}{3}t_2+\frac{15}{4}t_3+
\frac{24}{5}t_4 \\ + &
\frac{35}{6}t_5+\frac{48}{7}t_6+\frac{63}{8}t_7+\frac{80}{9}t_8+
\frac{99}{10}t_9+\frac{120}{11}t_{10}=e(S).
\end{split}
\end{equation*}
Noether formula gives $e(S)=12-K_S^2$, so we obtain
\begin{equation} \label{diofantina}
\begin{split}
& 41580K_S^2+41580t_1+73920t_2+103950t_3+133056t_4\\ + &
161700t_5+190080t_6+218295t_7+246400t_8\\+&
274428t_9+302400t_{10}=332640.
\end{split}
\end{equation}
We can check by direct computation that the only nonnegative
integers $K_S^2, \; t_1, \ldots, t_{10}$ which satisfy
(\ref{diofantina}) are
\begin{itemize}
\item $1 \leq K_S^2 \leq 8, \; \; t_1=8-K_S^2$;
\item $K_S^2=3, \; \; t_3=2$;
\item $K_S^2=2, \; \; t_1=1, \;  t_3=2$;
\item $K_S^2=1, \; \; t_1=2, \;  t_3=2$.
\end{itemize}
This completes the proof.
\end{proof}

\begin{proposition} \label{nocase7}
Let $\lambda \colon S \lr T=(C \times F)/G $ be as in Proposition \emph{\ref{poss-sing}}.
If $S$ is not a quasi-bundle, then $K_S^2 \leq 6$.
\end{proposition}
\begin{proof}
Since $S$ is not a quasi-bundle we have $K_S^2 \leq 7$. On the other
hand, if $K_S^2=7$ then $t_1=1$ and $t_n=0$ for $n \geq 2$. But this
is impossible by Remark \ref{singularities-2}.
\end{proof}
If $\chi(\mO_S)=1$ then Proposition \ref{ser:1} can be refined in
the following way.

\begin{proposition} \label{K2=6,5}
Let $S$ be as in Proposition \emph{\ref{poss-sing}} and let us assume $K_S^2=6$ or
$K_S^2=5$. Then
\begin{itemize}
\item $m_i\;$ divides $\;g(C)-1\;$ for all $i \in \{1, \ldots, r\}$, except at most one;
\item $n_j\;$ divides $\;g(F)-1\;$ for all $j \in \{1, \ldots, s \}$, except at most one.
\end{itemize}
\end{proposition}
\begin{proof}
Suppose $K_S^2=6$ or $K_S^2=5$. Then $T$ contains either $3$ or $2$
nodes (Proposition \ref{poss-sing}) and by Theorem \ref{Serrano} the
corresponding $(-2)-$curves must belong to the same fibre of $\beta
\colon S \lr F/G$. It follows that, for all $i$ except one, there is
a subgroup $H$ of $G$, isomorphic to $\mathbb{Z}_{m_i}$, which acts
freely on $C$. Now Riemann-Hurwitz formula applied to $C \lr C/H$
gives
\begin{equation*}
g(C)-1=m_i(g(C/H)-1),
\end{equation*}
 so $m_i$ divides $g(C)-1$. The second statement can be proven in
 the same way.
\end{proof}

Set $\mathbf{m}:=(m_1 , \ldots, m_r)$ and $\mathbf{n}:=(n_1 ,
\ldots, n_s)$, where we make the usual abbreviations such as $(2^3, 3^2)$.

\begin{proposition} \label{kappa5}
Let us assume $\chi(\mO_S)=1$ and $K_S^2=6$ or $K_S^2=5$. Then $g(F)=2$
 implies $\mathbf{n}=(2^1)$, whereas $g(C)=2$ implies
$\mathbf{m}=(2^1)$.
\end{proposition}
\begin{proof}
If $g(F)=2$ then Corollary \ref{ser:2} yields $\mathbf{n}=(2^s)$. On
the other hand, if $s \geq  2$ then Proposition $\ref{K2=6,5}$ implies
 that $2$  divides $g(F)-1=1$, a contradiction. An analogous proof works in the
 case $g(C)=2$.
\end{proof}

\section{Standard isotrivial fibrations with $p_g=q=1$. Building data} \label{structure theorem}
From now on we suppose that $\lambda \colon S \lr T=(C \times F)/G$ is
a standard isotrivial fibration with $p_g=q=1$,
such that $T$ has at worst RDPs.
Since $q=1$, we may assume that $E:=C/G$ is an
elliptic curve and that $F/G \cong \mathbb{P}^1$, that is $\mathfrak{g}'_1=0$ and
$\mathfrak{g}'_2=1$. Then the natural
projection $\alpha \colon S \lr E$ is the Albanese morphism of $S$
and $g_{\textrm{alb}}=g(F)$. Moreover by Proposition \ref{minimal} it follows
that $S$ is a minimal
model.
Let $\mathcal{V}=\{g_1, \ldots g_r \}$
be a generating vector for $G$ of type $(0\; | \; m_1, \ldots ,
m_r)$, inducing the $G-$cover $F \lr \mathbb{P}^1$ and let
$\mathcal{W}=\{\ell_1, \ldots \ell_s; \; h_1, h_2\}$ be a generating
vector of type $(1\; | \; n_1, \ldots , n_s)$ inducing $C \lr
E$. Then Riemann-Hurwitz formula implies
\begin{equation} \label{generi}
\begin{split}
2g(F)-2 &=|G|\bigg(-2+\sum_{i=1}^r \bigg( 1- \frac{1}{\;m_i} \bigg)
\bigg) \\
2g(C)-2 & =|G| \sum_{j=1}^s \bigg(1-\frac{1}{\;n_j} \bigg).
\end{split}
\end{equation}
\label{sec:prel}
\begin{proposition} \label{prel:prop1}
Let $\lambda \colon S \lr T=(C \times F)/G$ be a standard isotrivial
fibration with $p_g=q=1$, such that $T$ has at worst \emph{RDPs}.
Then
\begin{equation} \label{fundamental}
\frac{K_S^2}{4(g(F)-1)}= \Sn.
\end{equation}
\end{proposition}
\begin{proof}
Using Proposition \ref{invariantsS} and the second relation in \eqref{generi}  we
obtain
\begin{equation*}
\frac{|G| \cdot K_S^2}{4(g(F)-1)}=2(g(C)-1)=|G| \cdot \Sn,
\end{equation*}
so the claim follows.
\end{proof}

\begin{proposition} \label{nocase5}
The case $p_g=q=1$, $K_S^2=5$ does not occur.
\end{proposition}
\begin{proof}
If $K_S^2=5$ occurs, Proposition \ref{prel:prop1} gives
\begin{equation} \label{eq:1 no5}
(g(F)-1)\Sn=\frac{5}{4}.
\end{equation}
If $s \geq 2$ then $g(F)-1 \leq \frac{5}{4}$, hence $g(F)=2$. This
yields $\Sn=\frac{5}{4}$, hence $\mathbf{n}=(2^1, 4^1)$, which
contradicts Proposition \ref{kappa5}. Therefore we must have $s=1$,
i.e. $\mathbf{n}=(n^1)$. This implies
\begin{equation*}
\frac{5}{4}=(g(F)-1) \left(1-\frac{1}{n} \right) \geq
\frac{1}{2}(g(F)-1),
\end{equation*}
hence $g(F) \leq 3$. Using (\ref{eq:1 no5}), we obtain
$\left(1-\frac{1}{n} \right)=\frac{5}{4}$ if $g(F)=2$ and
$\left(1-\frac{1}{n} \right)=\frac{5}{8}$ if $g(F)=3$; but both
cases are impossible, because $n$ must be a positive integer.
\end{proof}

\begin{proposition} \label{nocase3}
The case $p_g=q=1$, $K_S^2=3$ does not occur.
\begin{proof}
If $K_S^2=3$ then either
 $g(F)=3$ or $g(F)=2$ (\cite{CaCi91}, \cite{CaCi93}). In the former case
 Proposition \ref{prel:prop1} implies
 $\Sn=\frac{3}{8}$, which is impossible. In the latter case we have
 $\Sn=\frac{3}{4}$, hence $\mathbf{n}=(4^1)$ which contradicts
 Corollary \ref{ser:2}.
\end{proof}
\end{proposition}

\begin{proposition} \label{cAse2}
If $p_g=q=1$ and $K_S^2=2$, then $T$ contains only nodes.
\end{proposition}
\begin{proof}
If $K_S^2=2$ we have $g(F)=2$ (\cite{Ca}, \cite{CaCi91},
\cite{CaCi93}). So the claim follows by Corollary \ref{nodal}.
\end{proof}

Summing up and using Proposition \ref{poss-sing} we obtain
\begin{proposition} \label{surfaces-are-nodal}
Let $\lambda \colon S \lr T=(C \times F)/G$ be a standard isotrivial
fibration with $p_g=q=1$, such that $T$ has at worst \emph{RDPs}.
Then $K_S^2$ is even and the only singularities of $T$ are $8-K_S^2$
nodes.
\end{proposition}
Now let us observe that the cyclic subgroups $\langle g_1 \rangle, \ldots ,\langle g_r
\rangle$ and their conjugates provide the non-trivial stabilizers of the action of $G$ on
$F$, whereas $\langle \ell_1 \rangle, \ldots, \langle \ell_s \rangle$ and their
conjugates provide the non-trivial stabilizers of the actions of $G$ on $C$. The
singularities of $T$ arise from the points in $C \times F$ with nontrivial stabilizer;
since the action of $G$ on $C \times F$ is the diagonal one, it follows that the set
 $\mathscr{S}$ of all nontrivial stabilizers for the action of $G$ on $C \times F$
 is given by
\begin{equation} \label{stabilizzatori}
\mathscr{S}= \bigg( \bigcup_{\sigma \in G} \bigcup_{i=1}^r \langle \sigma g_i \sigma^{-1}
\rangle \bigg) \cap \bigg( \bigcup_{\sigma \in G} \bigcup_{j=1}^s \langle \sigma \ell_j
\sigma^{-1} \rangle \bigg) \cap G^{\times}.
\end{equation}
Notice that  Proposition \ref{surfaces-are-nodal} implies that every element of
$\mathscr{S}$ has order $2$. Moreover the (reduced) fibre of the covering $C \times F \lr
T$ over each node has cardinality
 $\frac{|G|}{2}$, so the number of nodes of $T$ is given by
\begin{equation*}
8-K_S^2=t_1=\frac{2}{|G|}\sum_{h \in \mathscr{S}}|\textrm{Fix}_C(h)|
\cdot |\textrm{Fix}_F(h)|.
\end{equation*}
Proposition \ref{invariantsS} yields
\begin{equation} \label{Kappa2}
K_S^2=\frac{8(g(C)-1)(g(F)-1)}{|G|},
\end{equation}
so we can write down the basic equality
\begin{equation} \label{eq-fondam}
(g(C)-1)(g(F)-1)+\frac{1}{4}\sum_{h \in \mathscr{S}}|\textrm{Fix}_C(h)| \cdot
|\textrm{Fix}_F(h)|=|G|.
\end{equation}
We call $(G, \mathcal{V}, \mathcal{W})$ the \emph{building data} of
$S$. In fact, we have the following structure result.
\begin{proposition} \label{strukture}
Let $G$ be a finite group which is both $(0 \; |\;m_1, \ldots,
m_r)-$generated and \\$(1\; | \;n_1, \ldots, n_s)-$generated, with
generating vectors $\mathcal{V}=\{g_1, \ldots, g_r \}$ and
$\mathcal{W}=\{\ell_1, \ldots, \ell_s; \; h_1,h_2 \}$, respectively.
Denote by
\begin{equation*}
\begin{split}
f & \colon F \lr \mathbb{P}^1=F/G, \\   h & \colon C \lr E= C/G
\end{split}
\end{equation*}
the two $G-$coverings induced by $\mathcal{V}$ and $\mathcal{W}$ and
let $g(F), \; g(C)$ be the genera of $F$ and $C$, that are related
to $|G|, \; \mathbf{m}, \; \mathbf{n}$ by \emph{(\ref{generi})}. Finally,
 define $\mathscr{S}$ as in \emph{(\ref{stabilizzatori})}.
Assume moreover that
\begin{itemize}
\item $g(C) \geq 2, \;\; g(F) \geq2$;
 \item every element of $\mathscr{S}$ has order $2$;
 \item equality \emph{(\ref{eq-fondam})} is satisfied.
\end{itemize}
Then the quotient $T:=(C \times F)/G$ contains exactly $8-K_T^2$
nodes and its minimal desingularization $S$ is a minimal surface of
general type whose invariants are
\begin{equation*}
p_g(S)=q(S)=1, \quad  K_S^2=\frac{8(g(C)-1)(g(F)-1)}{|G|}.
\end{equation*}
 Conversely, every standard isotrivial fibration $S$, with $p_g(S)=q(S)=1$
 and such that $T$ has only \emph{RDPs}, arises in this way.
\end{proposition}
\begin{proof}
We have already shown that, if $\lambda \colon S \lr T=(C \times
F)/G$ is a standard isotrivial fibration with $p_g=q=1$, such that $T$ has
at worst RDPs, then the assumptions above must be satisfied. Vice
versa, if all the assumptions are satisfied then the quotient $T=(C
\times F)/G$ is a nodal surface with $q(T)=1$, whose number of nodes
is given by
\begin{equation*}
\begin{split}
t_1&=\frac{2}{|G|}\sum_{h \in \mathscr{S}}|\textrm{Fix}_C(h)| \cdot
|\textrm{Fix}_F(h)|  \\
&=\frac{2}{|G|}\cdot 4 \left(|G|-(g(C)-1)(g(F)-1)\right) \quad (\textrm{using}\; (\ref{eq-fondam}))  \\
&=8-\frac{8(g(C)-1)(g(F)-1)}{|G|}.
\end{split}
\end{equation*}
Let $S$ be the minimal desingularization  of $T$;
by using  Proposition \ref{invariantsS} and relation
\eqref{Kappa2} we obtain
\begin{equation*}
\begin{split}
e(S)&=\frac{1}{2}K_S^2+ \frac{3}{2}t_1 \\
&=\frac{1}{2}K_S^2+ \frac{3}{2}\left(8- K_S^2 \right) \\ &=12-K_S^2.
\end{split}
\end{equation*}
Thus Noether formula yields $\chi(\mO_S)=1$, that implies
$p_g(S)=q(S)=1$. Again by \eqref{Kappa2} we have $K_S^2 >0$, hence $S$ is
 a surface of general type, which must be minimal
by Proposition \ref{minimal}.
\end{proof}
\begin{remark}
The surface $S$ is a quasi bundle if and only if
$\mathscr{S}=\emptyset$ (see \cite[Proposition 7.2]{Pol08}).
\end{remark}

\section{Standard isotrivial fibrations with $p_g=q=1$. The abelian case}
\label{abelian case}
The aim of this section is to prove
\begin{theorem} \label{main-teo-ab}
Let $\lambda \colon S \lr T=(C \times F)/G$ be a standard
isotrivial fibration with $p_g=q=1$, which is not a
quasi-bundle, such that $T$ has only \emph{RDPs}. Assume in addition
that the group $G$ is abelian. Then $K_S^2=4,\; g(F)=2$ and there
are three cases:
\begin{itemize}
\item $g(C)=3, \quad G=\mZ_2 \times \mZ_2$
\item $g(C)=4, \quad G=\mZ_6$
\item $g(C)=7, \quad G=\mZ_2 \times \mZ_6$.
\end{itemize}
All possibilities occur.
\end{theorem}
The proof of Theorem \ref{main-teo-ab} will be a consequence of the following
results.

\begin{proposition} \label{ab:prop2}
If $G$ is abelian then
\begin{equation*}
K_S^2=4, \quad g(F)=2, \quad \mathbf{n}=(2^2).
\end{equation*}
\end{proposition}
\begin{proof}
Since $G$ is abelian, the $G-$cover $h \colon C \lr E$ is branched in
at least two points (Proposition \ref{ab-no1}); thus $\Sn \geq 1$.
By using Proposition \ref{prel:prop1} this gives
\begin{equation} \label{ab:eq1}
K_S^2 =4(g(F)-1) \Sn \geq 4(g(F)-1).
\end{equation}
Since $K_S^2 \leq 6$ (Proposition \ref{nocase7}), we obtain $g(F)=2$
and so $K_S^2 \geq 4$. Thus Proposition \ref{surfaces-are-nodal}
implies
 $K_S^2=6$ or $K_S^2=4$.
\begin{itemize}
\item If $K_S^2=6$ then $\Sn=\frac{3}{2}$, that is either $\mathbf{n}=(2^3)$ or
 $\mathbf{n}=(4^2)$; since $g(F)=2$, both possibilities
 contradict Proposition \ref{kappa5}.
\item If $K_S^2=4$ then $\Sn=1$, hence $\mathbf{n}=(2^2)$.
\end{itemize}
\end{proof}
\begin{corollary} \label{ab:cor1}
If $G$ is abelian then $|G|$ is even and $|G| \geq 4$.
\end{corollary}
\begin{proof}
By Propositions \ref{invariantsS} and \ref{ab:prop2} we obtain
 $|G|=2(g(C)-1)$, so $|G|$ is even. If $|G|=2$ then $g(C)=2$, so
 $S$ would be a minimal surface of general type with $p_g=q=1$, $K_S^2=4$
 and a rational pencil $|C|$ of genus $2$ curves; but this contradicts
 [Xi85, p.51]. Thus $|G| \geq 4$.
\end{proof}

\begin{proposition} \label{ab=4}
If $|G|=4$ then the only possibility is
\begin{equation*}
G=\mZ_2 \times \mZ_2, \quad \mathbf{m}=(2^5).
\end{equation*}
This case occurs.
\end{proposition}
\begin{proof}
If $|G|=4$ then Proposition \ref{ab:prop2} and relations
(\ref{generi}) imply $g(C)=3$ and
\begin{equation*}
-2+ \sum_{i=1}^r \left(1- \frac{1}{ \;m_i} \right)=\frac{1}{2},
\end{equation*}
so there are two possibilities:
\begin{itemize}
\item $\mathbf{m}=(2^2,4^2)$
\item $\mathbf{m}=(2^5)$.
\end{itemize}
First let us rule out the case $\mathbf{m}=(2^2,4^2)$. If it occurs,
 then $G=\mZ_4=\langle x \; | \; x^4=1 \rangle$. Up to automorphisms of $G$, we may assume
\begin{equation*}
\begin{split}
g_1&=g_2=x^2, \;\; g_3=x, \;\; g_4=x^3 \\
 \ell_1&= \ell_2=x^2.
\end{split}
\end{equation*}
Then $\mathscr{S}=\{ x^2\}$ and by using Corollary \ref{ordine2} we
obtain
\begin{equation*}
|\textrm{Fix}_F(x^2)|=6, \quad  |\textrm{Fix}_C(x^2)|=4.
\end{equation*}
It follows that equality (\ref{eq-fondam}) is not satisfied, so this case does not occur.
\\
It remains to show that the possibility $\mathbf{m}=(2^5)$ actually occurs.
 In this case $G=\mZ_2 \times \mZ_2$, because
 $\mathbb{Z}_4$ is not $(0\;|\; 2^5)-$generated. Our example is the
 following. \\ \\
$\bullet\;$ $G=\mZ_2 \times \mZ_2, \quad \mathbf{m}=(2^5), \quad g(C)=3$. \\
Set $\mathbb{Z}_2 \times \mZ_2= \langle x, y\; | \; x^2=y^2=[x,y]=1
\; \rangle$ and
\begin{equation*}
\begin{split}
g_1&=x, \; \; g_2=y, \; \;
g_3=g_4=g_5=xy\\
\ell_1&=\ell_2=x, \; \; h_1=h_2=y.
\end{split}
\end{equation*}
We have $\mathscr{S}= \{ x \}$ and by using Corollary \ref{ordine2}
we obtain
\begin{equation*}
|\textrm{Fix}_F(x)|=2, \quad  |\textrm{Fix}_C(x)|=4.
\end{equation*}
Equality (\ref{eq-fondam}) is satisfied, hence Proposition \ref{strukture}
implies that this case occurs.
\end{proof}

\begin{lemma} \label{cyclic}
If $G$ is cyclic then $m_1 \geq 3$.
\end{lemma}
\begin{proof}
If $G$ is cyclic then it contains a unique element $h$ of order $2$.
By Proposition \ref{ab:prop2} we have $\mathbf{n}=(2^2)$, hence $|\textrm{Fix}_C(h)|=2 \cdot
\frac{|G|}{2}=|G|$. On the other hand, if $m_1=2$ then
$|\textrm{Fix}_F(h)| \geq \frac{|G|}{2}$. Since $K_S^2=4$ we have
\begin{equation*}
4 =t_1=\frac{2}{|G|}\cdot |\textrm{Fix}_C(h)| \cdot
|\textrm{Fix}_F(h)| \geq |G|.
\end{equation*}
Thus $G=\mZ_4$, which contradicts Proposition \ref{ab=4}.
\end{proof}

\begin{proposition} \label{ab>4}
If $G$ is abelian and $|G| >4$ there are two possibilities:
\begin{itemize}
\item $G= \mZ_6, \quad \quad \quad \; \mathbf{m}=(3,6^2)$
\item $G=\mZ_2 \times \mZ_6, \quad \mathbf{m}=(2, 6^2)$.
\end{itemize}
Both cases occur.
\end{proposition}
\begin{proof}
The abelian group $G$ acts as a group
of automorphisms on the genus $2$ curve $F$ so that $F/G \cong \mathbb{P}^1$.
Let us look at Table \ref{2-abeliani} of  Appendix $A$. By using Corollary \ref{ab:cor1}
and Lemma \ref{cyclic} we may rule out
cases $(1a)$, $(1b)$, $(1c)$, $(1d)$, $(1e)$, $(1f)$, $(1h)$,
$(1i)$. It remains to show that cases $(1g)$ and $(1j)$  occur.
\\ \\
$\bullet$ Case $(1g)$. $G=\mZ_6, \quad \mathbf{m}=(3,6^2), \quad
g(C)=4$. \\
Set $\mZ_6=\langle x \; | \; x^6=1 \rangle$ and
\begin{equation*}
\begin{split}
g_1&=x^4, \;\; g_2=x, \;\; g_3=x \\
 \ell_1&= \ell_2=x^3, \;\; h_1=h_2=x.
\end{split}
\end{equation*}
Then $\mathscr{S}=\{ x^3 \}$ and
\begin{equation*}
|\textrm{Fix}_F(x^3)|=2, \quad  |\textrm{Fix}_C(x^3)|=6.
\end{equation*}
Equality (\ref{eq-fondam}) is satisfied, so this case occurs.\\\\
$\bullet$ Case $(1j)$. $G=\mZ_2 \times \mZ_6, \quad
\mathbf{m}=(2,6^2), \quad g(C)=7$. \\
Let $x$, $y$ be obvious generators of $G$ of order $2$
 and $6$, respectively, and set
\begin{equation*}
\begin{split}
g_1&=x, \; \; g_2= y^5, \; \;
g_3=xy \\
\ell_1&=\ell_2=y^3, \; \; h_1=x,  \; \; h_2=y.
\end{split}
\end{equation*}
Then $\mathscr{S}=\{y^3 \}$ and
\begin{equation*}
|\textrm{Fix}_F(y^3)|=2, \quad  |\textrm{Fix}_C(y^3)|=12.
\end{equation*}
Equality (\ref{eq-fondam}) is satisfied, so this
case
 occurs.\\ \\
The proof of Theorem \ref{main-teo-ab} is complete.
\end{proof}

\section{Standard isotrivial fibrations with $p_g=q=1$. The nonabelian case}
\label{general case}

By Proposition \ref{surfaces-are-nodal} we have $K_S^2=6$, $4$ or $2$. We deal with the
three cases separately.

\subsection{The case $K_S^2=6$} \label{K6}
\begin{proposition}
If $K_S^2=6$ then we have two possibilities:
\begin{itemize}
\item $g(F)=3, \quad \mathbf{n}=(4^1)$
\item $g(F)=4, \quad \mathbf{n}=(2^1)$.
\end{itemize}
\end{proposition}
\begin{proof}
Formula (\ref{fundamental}) in this case gives
\begin{equation} \label{eq:K6-1}
\frac{3}{2}=(g(F)-1) \Sn.
\end{equation}
If $s \geq 2$ then $3/2 \geq g(F)-1$, hence $g(F)=2$  which
contradicts Proposition \ref{kappa5}. Then $s=1$, i.e. $\mathbf{n}=(n^1)$.
Using (\ref{eq:K6-1}) we obtain $3/2 \geq 1/2(g(F)-1)$ which implies
$g(F) \leq 4$. The case $g(F)=2$ is impossible, otherwise
$1-1/n=3/2$; therefore either $g(F)=3$ or $g(F)=4$. Using again
(\ref{eq:K6-1}) we see that we have $\mathbf{n}=(4^1)$ in the
former case and $\mathbf{n}=(2^1)$ in the latter one.
\end{proof}

\begin{proposition} \label{complete-K6-g-3}
If $K_S^2=6$ and $g(F)=3$ there are precisely the following cases:
\begin{table}[H]
\begin{center}
\begin{tabular}{|c|c|c|}
\hline
$ $ & \verb|IdSmall| & $ $ \\
$G$ & \verb|Group|$(G)$ & $\mathbf{m}$ \\
\hline \hline
$\emph{SL}_2(\mathbb{F}_3)$ & $G(24,3)$ & $(3^2,6)$\\ \hline
$\mathbb{Z}_2 \ltimes (\mathbb{Z}_2 \times \mathbb{Z}_8)$ & $G(32,9)$ & $(2,4,8)$ \\ \hline
$\mathbb{Z}_2 \ltimes D_{2,8,5}$  & $G(32,11)$ & $(2,4,8)$ \\ \hline
$G(48,33)$ & $G(48,33)$ & $(2,3,12)$ \\ \hline
$\mathbb{Z}_3 \ltimes (\mathbb{Z}_4)^2$ & $G(48,3)$ &  $(3^2,4)$
\\ \hline
$\emph{PSL}_2(\mathbb{F}_7)$ & $G(168,42)$ & $(2,3,7)$ \\ \hline
\end{tabular}
\end{center}
\end{table}
\end{proposition}
\begin{proof}
By Proposition \ref{invariantsS} we have $3 \cdot |G|=8(g(C)-1)$, so $8$ divides
 $|G|$. The nonabelian group $G$ acts as a group of automorphisms on the genus $3$ curve
 $F$ so that $F/G \cong \mathbb{P}^1$.
In addition, since $\mathbf{n}=(4^1)$, it follows that $G$ must be
 $(1\; | \;4^1)-$generated.
Now let us look at Table \ref{3-nonabelian} of Appendix $A$; by using
Propositions
 \ref{Table-3-no-1-4-gen} and \ref{K2=6,5} we are only left with
 cases  $(3n)$, $(3q)$, $(3r)$, $(3t)$, $(3u)$, $(3w)$.  \\ \\
$\bullet$ Case $(3n)$. $G=\textrm{SL}_2(\mathbb{F}_3), \quad \mathbf{m}=(3^2,6), \quad g(C)=10$.\\
Set
\begin{equation*}
g_1=\left(
  \begin{array}{cc}
    1 & 0 \\
    1 & 1 \\
  \end{array}
\right) \quad
g_2=\left(
  \begin{array}{cc}
    1 & -1 \\
    0 &  \;\;1 \\
  \end{array}
\right)\quad
g_3=\left(
  \begin{array}{cc}
   \;\; 0 & 1 \\
   -1 & 1 \\
  \end{array}
\right)
\end{equation*}
\begin{equation*}
\;\; \ell_1=\left(
  \begin{array}{cc}
    \;\;1 & -1 \\
   -1 & -1 \\
  \end{array}
\right) \quad
 h_1=\left(
  \begin{array}{cc}
    1 & 0 \\
    1 & 1 \\
  \end{array}
\right)\quad
h_2=\left(
  \begin{array}{cc}
    1 & -1 \\
    0 & \;\;1 \\
  \end{array}
\right)
\end{equation*}
and $\ell= \left( \begin{array}{cc}
    -1 & \;\;0 \\
   \;\;0 & -1 \\
  \end{array}
\right)$. Since $(g_3)^3=(\ell_1)^2= \ell$ and $\ell \in Z(G)$ it follows
$\mathscr{S}=\textrm{Cl}(\ell)=\{\ell\}$. By using Corollary \ref{ordine2} we obtain
\begin{equation*}
|\textrm{Fix}_F(\ell)|=4, \quad |\textrm{Fix}_C(\ell)|=6
\end{equation*}
so equality (\ref{eq-fondam}) is satisfied and this case occurs.\\ \\
$\bullet$ Case $(3q)$. $G=\mathbb{Z}_2 \ltimes (\mathbb{Z}_2 \times
\mathbb{Z}_8)=G(32,9), \quad \mathbf{m}=(2,4,8), \quad g(C)=13$.\\
Set
\begin{equation*}
\begin{split}
g_1&=x, \; \; g_2=xz, \;\; g_3=z^7 \\
\ell_1&=yz^6, \; \; h_1=x, \;\; h_2=z.
\end{split}
\end{equation*}
Since $(\ell_1)^2=(g_3)^4=z^4$ and $z^4 \in Z(G)$,
we have $\mathscr{S}= \textrm{Cl}(z^4)=\{ z^4 \}$; moreover
$z^4 \notin \langle g_1 \rangle$ and  $z^4 \notin \langle g_2
\rangle$, so we obtain
\begin{equation*}
|\textrm{Fix}_F(z^4)|=4, \quad |\textrm{Fix}_C(z^4)|=8.
\end{equation*}
Thus equality (\ref{eq-fondam}) is satisfied and this case occurs.\\ \\
$\bullet$ Case $(3r)$. $G=\mathbb{Z}_2 \ltimes D_{2,8,5}=G(32,11),
\quad
\mathbf{m}=(2,4,8), \quad g(C)=13$.\\
Set
\begin{equation*}
\begin{split}
g_1&=x, \; \; g_2=xz, \;\; g_3=z^7 \\
\ell_1&=yz^6, \; \; h_1=x, \;\; h_2=z.
\end{split}
\end{equation*}
Since $(\ell_1)^2=(g_3)^4=z^4$ and $z^4 \in Z(G)$,
we have $\mathscr{S}= \textrm{Cl}(z^4)=\{ z^4 \}$; moreover $z^4
\notin \langle g_1 \rangle$ and  $z^4 \notin \langle g_2 \rangle$,
so we obtain
\begin{equation*}
|\textrm{Fix}_F(z^4)|=4, \quad |\textrm{Fix}_C(z^4)|=8.
\end{equation*}
Thus equality (\ref{eq-fondam}) is satisfied and this case occurs.\\ \\
$\bullet$ Case $(3t)$. $G=G(48,33), \quad \mathbf{m}=(2,3,12), \quad g(C)=19$.
\\
Set
\begin{equation*}
\begin{split}
g_1&=xz, \; \; g_2=zy^2, \;\; g_3=xy \\
\ell_1&=z, \; \; h_1=y^2, \;\; h_2=xy^2z.
\end{split}
\end{equation*}
Since $(\ell_1)^2=(g_3)^6=t$ and $t \in Z(G)$ we
have $\mathscr{S}= \textrm{Cl}(t)=\{t\}$; moreover $t \notin \langle
g_1 \rangle$, so we obtain
\begin{equation*}
|\textrm{Fix}_F(t)|=4, \quad |\textrm{Fix}_C(t)|=12.
\end{equation*}
Thus equality (\ref{eq-fondam}) is satisfied and this case occurs.\\ \\
$\bullet$ Case $(3u)$. $G=\mathbb{Z}_3 \ltimes (\mathbb{Z}_4)^2=G(48,3),
\quad \mathbf{m}=(3^2,4), \quad g(C)=19$. \\
Set
\begin{equation*}
\begin{split}
g_1&=x, \; \; g_2=x^{-1}y^{-1}, \;\; g_3=y \\
\ell_1&=yz^2, \; \; h_1=x, \;\; h_2=xyx.
\end{split}
\end{equation*}
We have $(\ell_1)^2=(g_3)^2=y^2$, so $\mathscr{S}=\textrm{Cl}(y^2)$.
One checks that $|C_G(y^2)|=16$, hence $|\mathscr{S}|=3$ (in fact,
$\mathscr{S}=\{y^2,\; xy^2x^2, \; x^2y^2x \}$). For every $h \in
\mathscr{S}$  we obtain
\begin{equation*}
|\textrm{Fix}_F(h)|=4, \quad |\textrm{Fix}_C(h)|=4.
\end{equation*}
Thus equality (\ref{eq-fondam}) is satisfied and this case occurs.\\ \\
$\bullet$ Case $(3w)$. $G=\textrm{PSL}_2(\mathbb{F}_7), \quad \mathbf{m}=(2,3,7), \quad g(C)=64$.\\
It is well known that $G$ can be embedded in
$S_8$; in fact \\ $G= \langle (375)(486), \; (126)(348) \rangle$.
Set
\begin{equation*}
\begin{split}
g_1&=(12)(34)(58)(67), \; \; g_2=(154)(367), \;\; g_3=(1247358) \\
\ell_1&=(1825)(3647), \; \; h_1=(2576348), \;\; h_2=(1673428).
\end{split}
\end{equation*}
The group $G$ contains $21$ elements of order $2$, which belong to a
unique conjugacy class (see \cite{CCPW} or \cite{Bar99}). Therefore
$\mathscr{S}=\textrm{Cl}(g_1)=\textrm{Cl}((\ell_1)^2)$ and
$|\mathscr{S}|=21$. It follows that for all $h \in \mathscr{S}$ we
have
\begin{equation*}
|\textrm{Fix}_F(h)|=4, \quad |\textrm{Fix}_C(h)|=2,
\end{equation*}
so equality (\ref{eq-fondam}) is satisfied and this case occurs.
Notice that in this example the Albanese fibre $F$ of $S$ is
isomorphic to the Klein plane quartic
$\{x_0x_1^3+x_1x_2^3+x_2x_0^3=0\} \subset \mathbb{P}^2$; in particular it is not
hyperelliptic. \\ \\
This completes the proof of Proposition \ref{complete-K6-g-3}.
\end{proof}

\begin{proposition} \label{complete-K6-g-4}
If $K_S^2=6$ and $g(F)=4$ there are precisely the following cases:
\begin{table}[ht!]
\begin{center}
\begin{tabular}{|c|c|c|}
\hline
$ $ & \verb|IdSmall| & $ $ \\
$G$ & \verb|Group|$(G)$ & $\mathbf{m}$ \\
\hline \hline $D_4$ & $G(8,3)$ & $(2^4,4)$\\ \hline $A_4$ &
$G(12,3)$ & $(2,3^3)$ \\ \hline
$D_{2,12,7} $ & $G(24,10)$ &
$(2,6,12)$ \\ \hline $\mZ_3 \times A_4$ & $G(36,11)$ &
$(3^2,6)$ \\ \hline $D_4 \ltimes (\mZ_3)^2$ & $G(72,40)$ & $(2,4,6)$ \\ \hline $S_5$ &
$G(120,34)$ & $(2,4,5)$ \\ \hline
\end{tabular}
\end{center}
\end{table}
\end{proposition}
\begin{proof}
By Proposition \ref{invariantsS} we have $|G|=4(g(C)-1)$, so $4$
divides
 $|G|$. Moreover, since $\mathbf{n}=(2^1)$, the group $G$ must be $(1\; | \;
2^1)-$generated. Now let us look at Table \ref{4-nonabelian} of Appendix
$A$; by using Propositions
 \ref{Table-4-no-1-2-gen} and \ref{K2=6,5} we are only left with
 cases  $(4c)$, $(4f)$, $(4r)$, $(4x)$, $(4ac)$, $(4ad)$. \\
 \\
$\bullet$ Case $(4c)$. $G=D_4, \quad \mathbf{m}=(2^4,4), \quad g(C)=3$. \\
Set
\begin{equation*}
\begin{split}
g_1&=x, \quad g_2=xy, \quad g_3=x, \quad g_4=xy^2, \quad g_5=y\\
\ell_1&=y^2, \quad h_1=y, \quad h_2=x.
\end{split}
\end{equation*}
We have $\mathscr{S}= \textrm{Cl}(\ell_1)=\{ y^2 \}$ and
\begin{equation*}
|\textrm{Fix}_F(y^2)|=2, \quad |\textrm{Fix}_C(y^2)|=4,
\end{equation*}
so equality (\ref{eq-fondam}) is satisfied and this case occurs. \\
\\
$\bullet$ Case $(4f)$. $G=A_4,\quad \mathbf{m}=(2, 3^3), \quad
g(C)=4$.
\\ Set
\begin{equation*}
\begin{split}
g_1&=(12)(34), \quad g_2=(134), \quad g_3=(134), \quad g_4=(123)\\
\ell_1&=(12)(34), \quad h_1=(123), \quad h_2=(124).
\end{split}
\end{equation*}
Then $\mathscr{S}= \textrm{Cl}(\ell_1)=\{(12)(34), \; (13)(24), \;
(14)(23) \}$. For all $h \in \mathscr{S}$ we have
\begin{equation*}
|\textrm{Fix}_F(h)|=2, \quad |\textrm{Fix}_C(h)|=2,
\end{equation*}
so equality (\ref{eq-fondam}) is satisfied and this case occurs. \\ \\
$\bullet$ Case $(4r)$. $G=D_{2,12,7}, \quad \mathbf{m}=(2,6,12), \quad g(C)=7$.\\
Set
\begin{equation*}
\begin{split}
g_1&=x, \quad g_2=y^5x, \quad g_3=y\\
\ell_1&=y^6, \quad h_1=x, \quad h_2=y.
\end{split}
\end{equation*}
We have $\ell_1=(g_3)^6$; since $\ell_1 \in Z(G)$ it follows
$\mathscr{S}=\textrm{Cl}(\ell_1)=\{ y^6 \}$. On the
other hand $\ell_1 \notin \langle g_1 \rangle$ and $\ell_1 \notin \langle g_2 \rangle$, so
we obtain
\begin{equation*}
|\textrm{Fix}_F(y^6)|=2, \quad |\textrm{Fix}_C(y^6)|=12.
\end{equation*}
Thus equality (\ref{eq-fondam}) is satisfied and this case occurs. \\ \\
$\bullet$ Case $(4x)$. $G=\mZ_3 \times A_4, \quad
\mathbf{m}=(3^2,6), \quad g(C)=10$.\\
Set $\mZ_3= \langle z \, | \, z^3=1 \rangle$ and
\begin{equation*}
\begin{split}
g_1&=(z, \, (123)), \quad g_2=(z, \, (234)),  \quad g_3=(z, \, (12)(34))\\
\ell_1&=(1, \, (12)(34)), \quad h_1=(1, (123)), \quad h_2=(z, \,(14)(23)).
\end{split}
\end{equation*}
Since $\ell_1=(g_3)^3$ we obtain $\mathscr{S}=\textrm{Cl}(\ell_1)$ and so $|\mathscr{S}|=3$.
For all $h \in \mathscr{S}$ we have
\begin{equation*}
|\textrm{Fix}_F(h)|=2, \quad |\textrm{Fix}_C(h)|=6,
\end{equation*}
so equality (\ref{eq-fondam}) is satisfied and this case occurs. \\ \\
$\bullet$ Case $(4ac)$. $G=D_4 \ltimes (\mZ_3)^2=G(72,40), \quad \mathbf{m}=(2,4,6),
\quad g(C)=19$.\\
Set
\begin{equation*}
\begin{split}
g_1&=xzy, \quad g_2=y, \quad g_3=y^2z^2x\\
\ell_1&=y^2, \quad h_1=xy, \quad h_2=xz.
\end{split}
\end{equation*}
We have $\ell_1=(g_2)^2$ and so  $\mathscr{S}=\textrm{Cl}(\ell_1)$;
since $|C_G(\ell_1)|=8$, it follows $|\mathscr{S}|=9$. Moreover
$g_1 \notin \textrm{Cl}(\ell_1)$ and $(g_3)^3 \notin \textrm{Cl}(\ell_1)$, hence
for all $h \in \mathscr{S}$ we
have
\begin{equation*}
|\textrm{Fix}_F(h)|=2, \quad |\textrm{Fix}_C(h)|=4.
\end{equation*}
Thus equality (\ref{eq-fondam}) is satisfied and this case occurs. \\ \\
$\bullet$ Case $(4ad)$. $G=S_5, \quad \mathbf{m}=(2,4,5), \quad g(C)=31$.\\
Set
\begin{equation*}
\begin{split}
g_1&=(12), \quad g_2=(1543), \quad g_3=(12345)\\
\ell_1&=(14)(35), \quad h_1=(145), \quad h_2=(1432).
\end{split}
\end{equation*}
We have $\ell_1=(g_2)^2$, hence $\mathscr{S}=\textrm{Cl}(\ell_1)$
and $|\mathscr{S}|=15$. For all $h \in \mathscr{S}$ we obtain
\begin{equation*}
|\textrm{Fix}_F(h)|=2, \quad |\textrm{Fix}_C(h)|=4,
\end{equation*}
so equality (\ref{eq-fondam}) is satisfied and this case occurs. \\ \\
This completes the proof of Proposition \ref{complete-K6-g-4}.
\end{proof}

\subsection{The case $K_S^2=4$}
\begin{proposition} \label{Kappa4}
If $K_S^2=4$ then we have two possibilities:
\begin{itemize}
\item $g(F)=2, \quad \mathbf{n}=(2^2)$
\item $g(F)=3, \quad \mathbf{n}=(2^1)$.
\end{itemize}
\end{proposition}
\begin{proof}
If $K_S^2=4$ then Proposition \ref{prel:prop1} gives
\begin{equation*}
(g(F)-1)\Sn=1.
\end{equation*}
If $s \geq 2$ then $g(F)-1 \leq 1$, which implies $g(F)=2$ and
$\mathbf{n}=(2^2)$. So we may assume $s=1$, i.e. $\mathbf{n}=(n^1)$.
In this case we have $\frac{1}{2}(g(F)-1) \leq 1$, then $g(F) \leq
3$. On the other hand, $g(F)=2$ gives $1-\frac{1}{n}=1$, a
contradiction; therefore $g(F)=3$ and $\mathbf{n}=(2^1)$.
\end{proof}
In Proposition \ref{ab:prop2} we have proven that if $G$ is abelian
then $K_S^2=4$ and $g_{\textrm{alb}}=2$. However
we can also obtain the same invariants
with nonabelian $G$:
\begin{proposition} \label{complete-2}
If $K_S^2=4, \; g(F)=2$ and $G$ is not abelian there are
precisely the following cases:
\begin{table}[ht!]
\begin{center}
\begin{tabular}{|c|c|c|}
\hline
$ $ & \verb|IdSmall| & $ $ \\
$G$ & \verb|Group|$(G)$ & $\mathbf{m}$ \\
\hline \hline
$S_3$  & $G(6,1)$ & $(2^2,3^2)$\\ \hline
$D_4$ & $G(8,3)$ &  $(2^3,4)$ \\ \hline
$D_6$  & $G(12,4)$ & $(2^3,3)$ \\ \hline
$D_{2,8,3}$ & $G(16,8)$ & $(2,4,8)$ \\ \hline
$\mathbb{Z}_2 \ltimes ((\mathbb{Z}_2)^2 \times
\mathbb{Z}_3)$ & $G(24,8)$ & $(2,4,6)$ \\ \hline
$\emph{GL}_2(\mathbb{F}_3)$ & $G(48,29)$ & $(2,3,8)$ \\ \hline
\end{tabular}
\end{center}
\end{table}
\end{proposition}
\begin{proof}
By Proposition \ref{invariantsS} we have $|G|=2(g(C)-1)$. Moreover,
 since $\mathbf{n}=(2^2)$,  it follows that $G$ is $(1\;
| \; 2^2)-$generated. Let us look at
Table \ref{2-nonabelian} of Appendix $A$. Using Proposition
\ref{no-generation-1-22}
we can rule out cases $(2b)$, $(2d)$ and $(2h)$.
Now we check the remaining possibilities. \\ \\
$\bullet$ Case ($2a$). $G=S_3, \quad \mathbf{m}=(2^2,3^2), \quad g(C)=4$.\\
Set
\begin{equation*}
\begin{split}
g_1&=(12), \; \; g_2=(12), \;\; g_3=(123), \; \;g_4=(132) \\
\ell_1&=\ell_2=(12), \; \; h_1=h_2=(13).
\end{split}
\end{equation*}
We have $\mathscr{S}=\textrm{Cl}(\ell_1)=\{(12), \; (13), \; (23)
\}$ and for every $h \in \mathscr{S}$ we obtain
\begin{equation*}
|\textrm{Fix}_F(h)|=2, \quad |\textrm{Fix}_C(h)|=2.
\end{equation*}
Thus equality (\ref{eq-fondam}) is satisfied and this case occurs.\\ \\
$\bullet$ Case ($2c$). $G=D_4, \quad \mathbf{m}=(2^3,4), \quad g(C)=5$. \\
Set
\begin{equation*}
\begin{split}
g_1&=x, \; \; g_2=xy, \;\; g_3=y^2, \; \;g_4=y \\
\ell_1&=\ell_2=x, \; \; h_1=h_2=y.
\end{split}
\end{equation*}
We have $\mathscr{S}=\textrm{Cl}(\ell_1)=\{x, xy^2\}$ and for every $h
\in \mathscr{S}$ we obtain
\begin{equation*}
|\textrm{Fix}_F(h)|=2, \quad |\textrm{Fix}_C(h)|=4.
\end{equation*}
Thus equality (\ref{eq-fondam}) is satisfied and this case occurs.\\ \\
$\bullet$ Case ($2e$). $G=D_6, \quad \mathbf{m}=(2^3,3), \quad g(C)=7$. \\
Set
\begin{equation*}
\begin{split}
g_1&=x, \; \; g_2=xy, \;\; g_3=y^3, \; \;g_4=y^2 \\
\ell_1&=xy, \;\;\ell_2=xy^5, \; \; h_1=x, \;\; h_2=y^2.
\end{split}
\end{equation*}
We have $\mathscr{S}=\textrm{Cl}(\ell_1)=\{xy, xy^3, xy^5\}$ and for
every $h \in \mathscr{S}$ we obtain
\begin{equation*}
|\textrm{Fix}_F(h)|=2, \quad |\textrm{Fix}_C(h)|=4.
\end{equation*}
Thus equality (\ref{eq-fondam}) is satisfied and this case occurs.\\ \\
$\bullet$ Case ($2f$). $G=D_{2,8,3}, \quad \mathbf{m}=(2,4,8), \quad g(C)=9$.\\
Set
\begin{equation*}
\begin{split}
g_1&=x, \; \; g_2=xy^7, \;\; g_3=y\\
\ell_1&=x, \;\;\ell_2=xy^6, \; \; h_1=x, \;\; h_2=y.
\end{split}
\end{equation*}
We have $\mathscr{S}=\textrm{Cl}(\ell_1)= \{x,\;xy^2,xy^4,xy^6 \}$.
Moreover $(g_2)^2=(g_3)^4=y^4$ and $y^4 \notin \mathscr{S}$, hence
for every $h \in \mathscr{S}$ we obtain
\begin{equation*}
|\textrm{Fix}_F(h)|=2, \quad |\textrm{Fix}_C(h)|=4.
\end{equation*}
Thus equality (\ref{eq-fondam}) is satisfied and this case occurs.\\ \\
$\bullet$ Case ($2g$). $G=\mathbb{Z}_2 \ltimes ((\mathbb{Z}_2)^2 \times \mathbb{Z}_3)=G(24,8), \quad \mathbf{m}=(2,4,6), \quad g(C)=13$.\\
Set
\begin{equation*}
\begin{split}
g_1&=x, \; \; g_2=wxz, \;\; g_3=zw\\
\ell_1&=\ell_2=x, \; \; h_1=z, \;\; h_2=w.
\end{split}
\end{equation*}
We have $\mathscr{S}=\textrm{Cl}(\ell_1)$; since
$C_G(\ell_1)=\langle x,y \rangle \cong \mZ_2 \times \mZ_2$, it
follows $|\mathscr{S}|=6$. Moreover $(g_2)^2 \notin
\textrm{Cl}(\ell_1)$ and $(g_3)^3 \notin \textrm{Cl}(\ell_1)$, so
for every $h \in \mathscr{S}$ we obtain
\begin{equation*}
|\textrm{Fix}_F(h)|=2, \quad |\textrm{Fix}_C(h)|=4.
\end{equation*}
Thus equality (\ref{eq-fondam}) is satisfied and this case occurs.\\ \\
$\bullet$ Case ($2i)$. $G=\textrm{GL}_2(\mathbb{F}_3), \quad \mathbf{m}=(2,3,8), \quad g(C)=25$.\\
Set
\begin{equation*}
\; \; \; \; \; g_1=\left(
  \begin{array}{cc}
    1 & \;\;1 \\
    0 & -1 \\
  \end{array}
\right) \quad
\; \; \; \; \;  g_2=\left(
  \begin{array}{cc}
    0 & -1 \\
    1 & -1 \\
  \end{array}
\right)\quad
\; \;   g_3=\left(
  \begin{array}{cc}
   -1 & \;\;1 \\
   -1 & -1 \\
  \end{array}
\right)
\end{equation*}
\begin{equation*}
\;\; \ell_1=\ell_2=\left(
  \begin{array}{cc}
    1 & \;\;1 \\
   0 & -1 \\
  \end{array}
\right) \quad
 h_1=\left(
  \begin{array}{cc}
    0 & -1 \\
    1 & -1 \\
  \end{array}
\right)\quad
h_2=\left(
  \begin{array}{cc}
    1 & 0 \\
    0 & 1 \\
  \end{array}
\right)
\end{equation*}
and $\ell= \left( \begin{array}{cc}
    -1 & \;\;0 \\
    \;\;0 & -1 \\
  \end{array}
\right)$. We have $\mathscr{S}=\textrm{Cl}(\ell_1)$ and
$C_G(\ell_1)\cong \mathbb{Z}_2 \times \mathbb{Z}_2$, hence
$|\mathscr{S}|=12$. Moreover $(g_3)^4=\ell \notin
\textrm{Cl}(\ell_1)$, so for all $h \in \mathscr{S}$ we obtain
\begin{equation*}
|\textrm{Fix}_F(h)|=2, \quad |\textrm{Fix}_C(h)|=4.
\end{equation*}
Thus equality (\ref{eq-fondam}) is satisfied and this case occurs. \\
\\ This completes the proof of Proposition \ref{complete-2}.
\end{proof}
 \newpage
\begin{proposition} \label{complete-K4-g-3}
If $K_S^2=4$ and $g(F)=3$ there are precisely the following cases:
\begin{table}[H]
\begin{center}
\begin{tabular}{|c|c|c|}
\hline
$ $ & \verb|IdSmall| & $ $ \\
$G$ & \verb|Group|$(G)$ & $\mathbf{m}$ \\
\hline \hline
$D_4$  & $G(8,3)$ & $(2^2,4^2)$\\ \hline
$D_4$ & $G(8,3)$ & $(2^5)$ \\ \hline
$A_4$ & $G(12,3)$ & $(2^2, 3^2)$ \\ \hline
$D_{2,8,5}$ & $G(16,6)$ & $(2,8^2)$ \\ \hline
$D_{4,4,-1}$ & $G(16,4)$ & $(4^3)$ \\ \hline
$\mathbb{Z}_2 \times A_4$ & $G(24,13)$ & $(2, 6^2)$ \\ \hline
\end{tabular}
\end{center}
\end{table}
\end{proposition}
\begin{proof}
By Proposition \ref{invariantsS} we have $|G|=4(g(C)-1)$, so $4$ divides
 $|G|$. Moreover, since $\mathbf{n}=(2^1)$, the group $G$ is $(1\; | \;
2^1)-$generated. Now let us look at Table \ref{3-nonabelian} of
 Appendix $A$; by using
Proposition \ref{Table-3-no-1-2-gen} we are only left with cases
$(3b)$, $(3c)$, $(3f)$, $(3g)$, $(3h)$, $(3m)$.  \\ \\
$\bullet$ Case $(3b)$. $G=D_4,  \; \; \mathbf{m}=(2^2,4^2), \; \; g(C)=3$. \\
Set
\begin{equation*}
\begin{split}
g_1&=x, \quad g_2=x, \quad g_3=y, \quad g_4=y^3  \\
\ell_1&=y^2, \quad h_1=x, \quad h_2=y.
\end{split}
\end{equation*}
We have $\mathscr{S}=\textrm{Cl}(\ell_1)=\{ y^2 \}$ and
\begin{equation*}
|\textrm{Fix}_F(y^2)|=4, \quad |\textrm{Fix}_C(y^2)|=4,
\end{equation*}
so equality (\ref{eq-fondam}) is satisfied and this case occurs. \\ \\
$\bullet$ Case $(3c)$. $G=D_4,  \; \; \mathbf{m}=(2^5), \; \; g(C)=3$. \\
Set
\begin{equation*}
\begin{split}
g_1&=y^2, \quad g_2=xy, \quad g_3=xy^3, \quad g_4=x, \quad g_5=x \\
\ell_1&=y^2, \quad h_1=x, \quad h_2=y.
\end{split}
\end{equation*}
We have $\mathscr{S}=\textrm{Cl}(\ell_1)=\{ y^2 \}$ and
\begin{equation*}
|\textrm{Fix}_F(y^2)|=4, \quad |\textrm{Fix}_C(y^2)|=4,
\end{equation*}
so equality (\ref{eq-fondam}) is satisfied and this case occurs. \\ \\
$\bullet$ Case $(3f)$. $G=A_4, \; \; \mathbf{m}=(2^2,3^2), \; \;
g(C)=4$. \\
Set
\begin{equation*}
\begin{split}
g_1&=(12)(34), \quad g_2=(12)(34), \quad g_3=(123), \quad g_4=(132) \\
\ell_1 &=(12)(34), \quad h_1=(123), \quad h_2=(124).
\end{split}
\end{equation*}
We have $\mathscr{S}=\textrm{Cl}(\ell_1)=\{(12)(34),\; (13)(24),\; (14)(23) \}$ and for all
$h \in \mathscr{S}$ we obtain
\begin{equation*}
|\textrm{Fix}_F(h)|=4, \quad |\textrm{Fix}_C(h)|=2,
\end{equation*}
so equality (\ref{eq-fondam}) is satisfied and this case occurs.\\ \\
$\bullet$ Case $(3g)$. $G=D_{2,8,5}, \; \; \mathbf{m}=(2,8^2), \; \;
g(C)=5$. \\
Set
\begin{equation*}
\begin{split}
g_1&=x, \quad g_2=xy^{-1}, \quad g_3=y \\
\ell_1&=y^4, \quad h_1=x, \quad h_2=y.
\end{split}
\end{equation*}
Since $\ell_1=(g_2)^4=(g_3)^4$ and $\ell_1 \in Z(G)$, it follows $\mathscr{S}=\textrm{Cl}(\ell_1)=\{y^4 \}$. We have
\begin{equation*}
|\textrm{Fix}_F(y^4)|=4, \quad |\textrm{Fix}_C(y^4)|=8,
\end{equation*}
so equality (\ref{eq-fondam}) is satisfied and this case occurs.\\ \\
$\bullet$ Case $(3h)$. $G=D_{4,4,-1}, \; \; \mathbf{m}=(4^3), \; \;
g(C)=5$ \\
Set
\begin{equation*}
\begin{split}
g_1&=x, \quad g_2=x^{-1}y^{-1}, \quad g_3=y \\
\ell_1&=y^2, \quad h_1=x, \quad h_2=y.
\end{split}
\end{equation*}
Since $\ell_1=(g_3)^2$ and $\ell_1 \in Z(G)$ we have
$\mathscr{S}=\textrm{Cl}(\ell_1)=\{ y^2\}$. Moreover $\ell_1 \notin
\langle g_1 \rangle$ and $\ell_2 \notin \langle g_2 \rangle$, so we
obtain
\begin{equation*}
|\textrm{Fix}_F(y^2)|=4, \quad |\textrm{Fix}_C(y^2)|=8.
\end{equation*}
Thus equality (\ref{eq-fondam}) is satisfied and this case occurs.\\ \\
$\bullet$ Case $(3m)$. $G=\mathbb{Z}_2 \times A_4, \; \;
\mathbf{m}=(2,6^2), \; \; g(C)=7$. \\
Let $\mathbb{Z}_2= \langle z \; | \; z^2=1 \rangle$ and set
\begin{equation*}
\begin{split}
g_1&=(1, \; (12)(34)), \quad g_2=(z, \;(123)), \quad g_3=(z, \;(234)) \\
\ell_1 &=(1, \;(12)(34)), \quad h_1=(z, \; (123)), \quad h_2=(z,
\;(124)).
\end{split}
\end{equation*}
We have $\mathscr{S}=\textrm{Cl}(\ell_1)=\{(1,\;(12)(34)),\; (1,
\;(13)(24)),\; (1,\;(14)(23)) \}$. For all $h \in \mathscr{S}$ we
obtain
\begin{equation*}
|\textrm{Fix}_F(h)|=4, \quad |\textrm{Fix}_C(h)|=4,
\end{equation*}
so equality (\ref{eq-fondam}) is satisfied and this case occurs. \\
\\
This completes the proof of Proposition \ref{complete-K4-g-3}.
\end{proof}

\subsection{The case $K_S^2=2$}

\begin{lemma} \label{K2-2^1}
If $K_S^2=2$ then $\mathbf{n}=(2^1)$.
\end{lemma}
\begin{proof}
If $K_S^2=2$ we have $g(F)=2$ (\cite{Ca}, \cite{CaCi91},
\cite{CaCi93}). Therefore by
 Proposition \ref{prel:prop1} we obtain $\Sn=\frac{1}{2}$, that is
 $\mathbf{n}=(2^1)$.
\end{proof}

\begin{proposition} \label{completeK2}
If $K_S^2=2$ there are precisely the following possibilities:
\begin{table}[ht!]
\begin{center}
\begin{tabular}{|c|c|c|}
\hline
$ $ & \verb|IdSmall| & $ $ \\
$G$ & \verb|Group|$(G)$ & $\mathbf{m}$ \\
\hline \hline
$Q_8$  & $G(8,4)$ &  $(4^3)$\\
\hline
$D_4$ & $G(8,3)$ & $(2^3,4)$ \\
\hline
\end{tabular}
\end{center}
\end{table}
\end{proposition}
\begin{proof}
Proposition \ref{invariantsS} yields $|G|=4(g(C)-1)$, so  $4$
divides $|G|$. Since $\mathbf{n}=(2^1)$,
$G$ must be $(1 \; | \; 2^1)-$generated.
Now let us look at Table \ref{2-nonabelian} of Appendix $A$.
By using Proposition \ref{no-generation-1-21}
 we may rule out cases $(2d)$, $(2e)$, $(2f)$, $(2g)$, $(2h)$ and $(2i)$,  so the proof will be
complete if we show that cases $(2b)$ and $(2c)$  occur.
\\ \\
$\bullet$ Case $(2b)$. $G=Q_8,\quad \mathbf{m}=(4^3), \quad g(C)=3$. \\ Set
\begin{equation*}
\begin{split}
g_1&=j, \; \; g_2=i, \;\; g_3=k \\
\ell_1&=-1, \; \; h_1=i, \; \; h_2=j.
\end{split}
\end{equation*}
We have $\mathscr{S}= \textrm{Cl}(\ell_1)=\{-1 \}$ and
\begin{equation*}
|\textrm{Fix}_F(-1)|=6, \quad |\textrm{Fix}_C(-1)|=4,
\end{equation*}
so equality (\ref{eq-fondam}) is satisfied and this case occurs.\\ \\
$\bullet$ Case $(2c)$. $G=D_4, \quad \mathbf{m}=(2^3,4), \quad g(C)=3$. \\
Set
\begin{equation*}
\begin{split}
g_1&=xy^2, \; \; g_2=xy^3, \;\; g_3=y^2, \; \;g_4=y \\
\ell_1&=y^2, \; \; h_1=x, \; \; h_2=y.
\end{split}
\end{equation*}
We have $\mathscr{S}=\textrm{Cl}(y^2)=\{ y^2 \}$ and
\begin{equation*}
|\textrm{Fix}_F(y^2)|=6, \quad |\textrm{Fix}_C(y^2)|=4,
\end{equation*}
so equality (\ref{eq-fondam}) is satisfied and this case occurs.
\end{proof}

Proposition \ref{completeK2} shows that there exist two families of
 standard isotrivial fibrations with $p_g=q=1, \; K_S^2=2$. The
first family, that we denote by $\mathfrak{M}_{D_4}$, has dimension $2$
because it depends on the choice of four points on $\mathbb{P}^1$ and one point
on $E$ (up to projective equivalence);  the second family, that we
denote by $\mathfrak{M}_{Q_8}$, has dimension $1$ because it depends on the choice
of three points on $\mathbb{P}^1$ and one point on $E$.
Now we can provide a geometric description of $\mathfrak{M}_{D_4}$ and
 $\mathfrak{M}_{Q_8}$; to this purpose, let us recall some facts about
surfaces of general type with $p_g=q=1, \: K_S^2=2$ (see \cite{Ca}
and \cite{CaCi91} for further details). Let $(E, \oplus, 0)$ be
an elliptic curve $E$  with group law $\oplus$ and identity element
$0$, and let
\begin{equation*}
E^{(2)}=\textrm{Sym}^2(E)=\{x+y \; | \; x,y \in E \}
\end{equation*}
be its double symmetric product. Then the Abel-Jacobi map $E^{(2)}
\lr E$, $x+y \lr x \oplus y$  gives to $E^{(2)}$  the structure of a
$\mathbb{P}^1-$bundle over $E$. For any $a \in E$, let us consider the following
divisors on $E^{(2)}$:
\begin{equation*}
\begin{split}
\mathfrak{f}_a&:=\{x+y \in E \;|\; x \oplus y =a \}; \\
\mathfrak{h}_a&:= \{x+a \;|\; x \in E \}.
\end{split}
\end{equation*}
In both cases the corresponding
algebraic equivalence classes do not depend on $a$, hence we
may denote them by $\mathfrak{f}$ and $\mathfrak{h}$, respectively.
We have $\textrm{NS}(E^{(2)})=\mathbb{Z}\; \mathfrak{f} \oplus
\mathbb{Z}\; \mathfrak{h}$. The antibicanonical system
$|-2K_{E^{(2)}}|=|4 \mathfrak{h}_0 -2 \mathfrak{f}_0|$ is a linear pencil,
whose general elements are smooth
elliptic curves of the form
\begin{equation*}
\mathfrak{b}_a:= \{x+(x \oplus a)\; | \; x \in E\}, \quad a \oplus a \neq 0.
\end{equation*}
If $\div a$ denotes the inverse element of $a \in E$, we have
$\mathfrak{b}_a=\mathfrak{b}_{\div a}$.
It follows that the  singular members of $|-2K_{E^{(2)}}|$ are
precisely  the three double curves
$2\mathfrak{b}_{\xi_1}, \; 2\mathfrak{b}_{\xi_2}, \;
2\mathfrak{b}_{\xi_3}$, where the $\xi_i$ are the three
$2-$torsion points of $E$ different from $0$. The  $\mathfrak{b}_{\xi_i}$
are three divisors on $E^{(2)}$ which are algebraically but not
linearly equivalent to $2 \mathfrak{h}_0 - \mathfrak{f}_0$ (in fact,
$\mathfrak{b}_{\xi_i} \in |2 \mathfrak{h}_0 -
\mathfrak{f}_{\xi_i}|$). In \cite{Ca} it is shown that any surface
$S$ of general type with $p_g=q=1, \; K_S^2=2$ is a double cover of
$E^{(2)}$ branched along a divisor $\mathfrak{B}$
 algebraically equivalent to $6 \mathfrak{h} -2 \mathfrak{f}$ and having
at worst simple singularities.
 In particular the Albanese pencil $\{F \}$ of $S$ is the pullback of the
 ruling $\{\mathfrak{f}\}$ of $E^{(2)}$. Since the group of translations
 of $E$ acts transitively on the set of linear equivalence classes of divisors
 algebraically equivalent to $6 \mathfrak{h} -2 \mathfrak{f}$, we may assume
$\mathfrak{B} \in |6 \mathfrak{h}_0 -2 \mathfrak{f}_0|$. Therefore the
surfaces in $\mathfrak{M}_{D_4}$ and $\mathfrak{M}_{Q_8}$ must
 correspond to special curves with six nodes in the linear
 system $|6 \mathfrak{h}_0 -2 \mathfrak{f}_0|$. Indeed we can prove

\begin{proposition} \label{descriptionK2=2}
Let $S$ be the double cover of $E^{(2)}$
branched along a curve $\mathfrak{B} \in |6 \mathfrak{h}_0 -2 \mathfrak{f}_0|$.
Then the following holds.
\begin{itemize}
\item[$(i)$] If $S \in  \mathfrak{M}_{D_4}$ we have
\begin{equation*}
\mathfrak{B}=\mathfrak{B}'+\mathfrak{b}_{\xi_i}+\mathfrak{f}_{\xi_i},
\end{equation*}
where $\mathfrak{B}' \in  |-2K_{E^{(2)}}|$ and $i \in \{1,2,3\}$.
\item[$(ii)$] If $S \in \mathfrak{M}_{Q_8}$ we have
\begin{equation*}
\mathfrak{B}=\mathfrak{b}_{\xi_1}+ \mathfrak{b}_{\xi_2} + \mathfrak{b}_{\xi_3}   +\mathfrak{f}_0.
\end{equation*}
\end{itemize}
In both cases the isotrivial fibration $|C|$ of $S$ is obtained as
the pullback of the antibicanonical
pencil $|-2K_{E^{(2)}}|$ of $E^{(2)}$.
\end{proposition}
\begin{proof}
Since $C^2=0$ and $CF=8$, it follows that the image of $|C|$
 in $E^{(2)}$ via the double cover  $S\lr
E^{(2)}$ is a linear pencil whose general element $\mathfrak{c}$
verifies $\mathfrak{c}^2=0, \; \mathfrak{c} \mathfrak{f}=4$.
This implies $|\mathfrak{c}|=|-2K_{E^{(2)}}|$
(\cite{CaCi93}, p.404). Moreover, since $\mathbf{n}=(2^1)$, exactly one component
of $\mathfrak{B}$ is algebraically equivalent to $\mathfrak{f}$. If $S
\in \mathfrak{M}_{D_4}$ then $\mathbf{m}=(2^2,4)$, so exactly one of the curves
$\mathfrak{b}_{\xi_i}$ is contained in $\mathfrak{B}$; this implies $(i)$.
If $S \in \mathfrak{M}_{Q_8}$ then $\mathbf{m}=(4^3)$, so all the
 $\mathfrak{b}_{\xi_i}$ are components of $\mathfrak{B}$; this implies
 $(ii)$.
 \end{proof}
Notice that in both cases all the components of $\mathfrak{B}$ not
 contained in $\{ \mathfrak{f} \}$ are
invariant under translation in $E^{(2)}$; this explains why the Albanese
pencil of $S$ turns out to be isotrivial.
\newpage
\section*{Appendix $A$} \label{appendix A}
This appendix contains the classification of finite groups of automorphisms
 acting on Riemann surfaces of genus $2$, $3$ and $4$ so that the quotient is
 isomorphic to $\mathbb{P}^1$.
In the last two cases we listed only the nonabelian groups.
Tables \ref{2-abeliani}, \ref{2-nonabelian} and
 \ref{3-nonabelian} are adapted from [Br90, pages 252, 254, 255], whereas
Table \ref{4-nonabelian} is adapted from [Ki03, Theorem 1] and \cite{Vin00}.
For every $G$ we give a presentation, the vector $\mathbf{m}$ of branching data  and
the \verb|IdSmallGroup|$(G)$, that is the number of $G$ in the GAP4 database of small groups.
The author  wishes to thank
S. A. Broughton who kindly communicated to him that the group $G(48,33)$
(Table \ref{3-nonabelian}, case $(3t)$)
was missing in \cite{Br90}.
\begin{table}[ht!]
\begin{center}
\begin{tabular}{|c|c|c|c|}
\hline
$ $ & $ $ & \verb|IdSmall| & $ $ \\
$\textrm{Case}$ & $G$ & \verb|Group|$(G)$ &  $\mathbf{m}$ \\
\hline \hline
$(1a)$ & $\mZ_2$ & $G(2,1)$ & $(2^6)$ \\
\hline
$(1b)$ & $\mZ_3$ & $G(3,1)$ & $(3^4)$ \\
\hline
$(1c)$ & $\mZ_4$ & $G(4,1)$ & $(2^2, 4^2)$ \\
\hline
$(1d)$ & $\mZ_2 \times \mZ_2$ & $G(4,2)$ & $(2^5)$ \\
\hline
$(1e)$  & $\mZ_5$  & $G(5,1)$ & $(5^3)$\\
\hline
 $(1f)$  & $\mZ_{6}$ & $G(6,2)$ & $(2^2,3^2)$ \\
\hline
$(1g)$  & $\mZ_6$ & $G(6,2)$ & $(3,6^2)$ \\
\hline
$(1h)$  & $\mZ_8$ & $G(8,1)$ & $(2,8^2)$ \\
\hline
$(1i)$  & $\mZ_{10}$ & $G(10,2)$ & $(2,5,10)$ \\
\hline
$(1j)$  & $\mZ_2 \times \mZ_6$ & $G(12,5)$ & $(2,6^2)$ \\
\hline
\end{tabular}
\end{center}
\caption{Abelian groups of automorphisms acting with rational
quotient
 on Riemann
  surfaces of genus $2$}
\label{2-abeliani}
\end{table}

\begin{table}[ht!]
\begin{center}
\begin{tabular}{|c|c|c|c|c|}
\hline
$ $ & $ $ &  \verb|IdSmall| & $ $ & $ $ \\
${\textrm{Case}}$ & $G$ & \verb|Group|$(G)$ & $\mathbf{m}$ & $ \textrm{Presentation}$ \\
\hline \hline
 $(2a)$  & $S_3$  & $G(6,1)$ & $(2^2,3^2)$ & $\langle x,y\;|\;x=(123), \; y=(12)
 \rangle$  \\
\hline
$$ & $$ & $$ & $$  & $\langle i,j,k\; | \; i^2=j^2=k^2=-1$,
\\ $(2b)$ & $Q_8$ & $G(8,4)$ & $(4^3)$ & $ij=k, \; jk=i, \; ki=j \rangle$ \\
\hline
$(2c)$ & $D_4$ & $G(8,3)$ & $(2^3,4)$ & $\langle  x,y\; | \; x^2=y^4=1,
 \; xyx^{-1}=y^{-1} \rangle$ \\
\hline
$(2d)$ & $D_{4,3,-1}$ & $G(12,1)$ & $(3,4^2)$ & $\langle x,y\; | \; x^4=y^3=1,
 \; xyx^{-1}=y^{-1} \rangle$ \\
\hline
$(2e)$ & $D_6$ & $G(12,4)$ & $(2^3,3)$ & $\langle  x,y \; | x^2=y^6=1,
 \; xyx^{-1}=y^{-1} \rangle$ \\
\hline
$(2f)$ & $D_{2,8,3}$ & $G(16,8)$ & $(2,4,8)$ & $\langle  x,y \; | \;x^2=y^8=1,
 \; xyx^{-1}=y^3 \rangle$ \\
\hline
$$ & $$ & $$ & $$ & $\langle x,y,z,w \; | \; x^2=y^2=z^2=w^3=1,$
\\
$(2g)$ & $G=\mathbb{Z}_2 \ltimes ((\mathbb{Z}_2)^2 \times
\mathbb{Z}_3)$ & $G(24,8)$ & $(2,4,6)$ & $[y,z]=[y,w]=[z,w]=1,$
\\
$$ & $$ & $$ & $$ & $xyx^{-1}=y,\; xzx^{-1}=zy, \; xwx^{-1}=w^{-1} \rangle$ \\
\hline
$(2h)$ &  $\textrm{SL}_2(\mathbb{F}_3)$ & $G(24,3)$ & $(3^2,4)$ & $ \langle x,y \; | \;
 x=
\left(
              \begin{array}{cc}
                1 & 1 \\
                0 & 1 \\
              \end{array}
            \right), \;
y=
\left(
              \begin{array}{cc}
                \;\;0 & \;\;1 \\
                -1 & -1 \\
              \end{array}
            \right)
\rangle$ \\
\hline
$(2i)$ &  $\textrm{GL}_2(\mathbb{F}_3)$ & $G(48,29)$ & $(2,3,8)$ & $ \langle x,y \; | \;
 x=
\left(
              \begin{array}{cc}
                1 & \;\;1 \\
                0 & -1 \\
              \end{array}
            \right), \;
 y=
\left(
              \begin{array}{cc}
                0 & -1 \\
                1 & -1 \\
              \end{array}
            \right)
\rangle$ \\
\hline
\end{tabular}
\end{center}
\caption{Nonabelian groups of automorphisms acting with rational
quotient on Riemann surfaces of genus $2$.} \label{2-nonabelian}
\end{table}

\newpage

\begin{table}[ht!]
\begin{center}
\begin{tabular}{|c|c|c|c|c|}
\hline
$ $ & $ $ &  \verb|IdSmall| & $ $ & $ $ \\
${\textrm{Case}}$ & $G$ & \verb|Group|$(G)$ & $\mathbf{m}$ & $ \textrm{Presentation}$ \\
\hline \hline $(3a)$  & $S_3$  & $G(6,1)$ & $(2^4,3)$ & $\langle x,y\;|\;x=(12), \;
y=(123)
 \rangle$  \\
 \hline
$(3b)$ & $D_4$ & $G(8,3)$ & $(2^2,4^2)$ & $\langle  x,y\; | \; x^2=y^4=1,
 \; xyx^{-1}=y^{-1} \rangle$ \\
\hline $(3c)$ & $D_4$ & $G(8,3)$ & $(2^5)$ & $\langle  x,y\; | \; x^2=y^4=1,
 \; xyx^{-1}=y^{-1} \rangle $ \\
\hline
 $(3d)$ & $D_{4,3,-1}$ & $G(12,1)$ & $(4^2,6)$ & $\langle x,y\; | \; x^4=y^3=1,
 \; xyx^{-1}=y^{-1} \rangle$ \\
\hline $(3e)$ & $D_6$ & $G(12,4)$ & $(2^3,6)$ & $\langle  x,y\; | \; x^2=y^6=1,
 \; xyx^{-1}=y^{-1} \rangle$ \\
 \hline
$(3f)$ & $A_4$ & $G(12,3)$ & $(2^2, 3^2)$ & $\langle x,y \; | \;
x=(12)(34), \; y=(123) \rangle$ \\
\hline $(3g)$ & $D_{2,8,5}$ & $G(16,6)$ & $(2,8^2)$ & $\langle x,y\; | \; x^2=y^8=1,
 \; xyx^{-1}=y^5 \rangle$ \\
\hline $(3h)$ & $D_{4,4,-1}$ & $G(16,4)$ & $(4^3)$ & $\langle x,y\; | \; x^4=y^4=1,
 \; xyx^{-1}=y^{-1} \rangle$ \\
\hline $(3i)$ & $\mathbb{Z}_2 \times D_4$ & $G(16,11)$ & $(2^3,4)$ & $\langle z \; | \;
z^2=1 \rangle \times   \langle  x,y\; | \; x^2=y^4=1,
 \; xyx^{-1}=y^{-1} \rangle $ \\
\hline $ $ & $ $ & $ $ & $ $ & $\langle x,y,z \; | \; x^2=y^2=z^4=1,$\\
$(3j)$ & $\mathbb{Z}_2 \ltimes (\mathbb{Z}_2 \times
\mathbb{Z}_4) $ & $G(16,13)$ & $(2^3,4)$ & $[x,z]=[y,z]=1, \; xyx^{-1}=yz^2 \rangle$ \\
\hline $(3k)$ & $D_{3,7,2}$ & $G(21,1)$ & $(3^2,7)$ & $\langle x,y\; | \; x^3=y^7=1,
 \; xyx^{-1}=y^2 \rangle$ \\
\hline $(3l)$ & $D_{2,12,5}$ & $G(24,5)$ & $(2,4,12)$ & $\langle x,y\; | \; x^2=y^{12}=1,
 \; xyx^{-1}=y^5 \rangle$ \\
\hline $(3m)$ & $\mathbb{Z}_2 \times A_4$ & $G(24,13)$ & $(2,6^2)$ & $ \langle z \; | \;
z^2=1 \rangle \times \langle x,y \; | \; x=(12)(34), \; y=(123) \rangle
$ \\
\hline $(3n)$ & $\textrm{SL}_2(\mathbb{F}_3)$ & $G(24,3)$ & $(3^2,6)$ & $ \langle x,y \;
| \;
 x=
\left(
              \begin{array}{cc}
                1 & 1 \\
                0 & 1 \\
              \end{array}
            \right), \;
y= \left(
              \begin{array}{cc}
                \;\;0 & \;\;1 \\
                -1 & -1 \\
              \end{array}
            \right)
\rangle$ \\
\hline
$(3o)$ & $S_4$ & $G(24,12)$ & $(3,4^2)$ & $\langle x,y \; | \; x=(1234), \; y=(12)   \rangle$ \\
\hline
$(3p)$ & $S_4$ & $G(24,12)$ & $(2^3,3)$ & $\langle x,y \; | \; x=(1234), \; y=(12)   \rangle$ \\
\hline $ $ & $ $ & $ $ & $ $ & $\langle x,y,z \; | \; x^2=y^2=z^8=1,$\\
$(3q)$ & $\mathbb{Z}_2 \ltimes(\mathbb{Z}_2 \times
\mathbb{Z}_8)$ & $G(32,9)$ & $(2,4,8)$ & $[x,y]=[y,z]=1, \; xzx^{-1}=yz^3 \rangle$ \\
\hline
$ $ & $ $ & $ $ & $ $ & $\langle x,y,z \; | \; x^2=y^2=z^8=1,$\\
$(3r)$ & $\mathbb{Z}_2 \ltimes D_{2,8,5}$ & $G(32,11)$ & $(2,4,8)$ & $yzy^{-1}=z^5, \; xyx^{-1}=yz^4, \; xzx^{-1}=yz^3 \rangle$ \\
\hline $(3s)$ & $\mathbb{Z}_2 \times S_4$ & $G(48,48)$ & $(2,4,6)$ & $ \langle z \; | \;
z^2=1 \rangle \times \langle x,y \; | \; x=(12), \; y=(1234) \rangle
$ \\
\hline $ $ & $ $ & $ $ & $ $ & $ \langle x,y,z,w,t \; | \; x^2=z^2=w^2=t,
\; y^3=1,\; t^2=1, $ \\
$(3t)$ & $G(48,33)$ & $G(48,33)$ & $(2,3,12)$ & $yzy^{-1}=w, \; ywy^{-1}=zw, \; zwz^{-1}=wt$,\\
$ $ & $ $ & $ $ & $ $ & $[x,y]=[x,z]=1 \rangle$ \\   \hline
$ $ & $ $ & $ $ & $ $ & $\langle x,y,z \; | \; x^3=y^4=z^4=1,$\\
$(3u)$ & $\mathbb{Z}_3 \ltimes(\mathbb{Z}_4)^2$ & $G(48,3)$ &
$(3^2,4)$ & $[y,z]=1, \; xyx^{-1}=z, \; xzx^{-1}=(yz)^{-1} \rangle$ \\
\hline
$$ & $$ & $$ & $$ & $\langle x,y,z,w \; | \; x^2=y^3=z^4=w^4=1,$
\\
$(3v)$ & $S_3 \ltimes (\mathbb{Z}_4)^2$ & $G(96,64)$ & $(2,3,8)$ & $[z,w]=1, \;
xyx^{-1}=y^{-1}, \; xzx^{-1}=w,$
\\
$$ & $$ & $$ & $$ & $xwx^{-1}=z,\; yzy^{-1}=w, \; ywy^{-1}=(zw)^{-1} \rangle$ \\\hline
$(3w)$ & $\textrm{PSL}_2(\mathbb{F}_7)$ & $G(168,42)$ & $(2,3,7)$ &
$\langle x,y \; | \; x=(375)(486), \; y=(126)(348)\rangle$ \\
\hline
\end{tabular}
\end{center}
\caption{Nonabelian groups of automorphisms acting with rational
quotient
  on Riemann surfaces of genus $3$.}
\label{3-nonabelian}
\end{table}

\newpage

\begin{table}[ht!]
\begin{center}
\begin{tabular}{|c|c|c|c|c|}
\hline
$ $ & $ $ &  \verb|IdSmall| & $ $ & $ $ \\
${\textrm{Case}}$ & $G$ & \verb|Group|$(G)$ & $\mathbf{m}$ & $ \textrm{Presentation}$ \\
\hline \hline $(4a)$ & $S_3$ & $G(6,1)$ & $(2^6)$ & $\langle x,y\;|\;x=(12), \;
y=(123) \rangle$  \\
\hline $(4b)$ & $S_3$ & $G(6,1)$ & $(2^2,3^3)$ & $\langle x,y\;|\;x=(12), \;
y=(123) \rangle$  \\
\hline $(4c)$ & $D_4$ & $G(8,3)$ & $(2^4,4)$ & $\langle  x,y\; | \; x^2=y^4=1,
 \; xyx^{-1}=y^{-1} \rangle$ \\
\hline
$ $ & $ $ & $ $ & $ $ &$ \langle i,j,k,-1\; | \; i^2=j^2=k^2=-1,$\\
$(4d)$ & $Q_8$ & $G(8,4)$ & $(2,4^3)$ & $ij=k, \; jk=i, \; ki=j \rangle$ \\
\hline $(4e)$ & $D_5$ & $G(10,1)$ & $(2^2,5^2)$ & $\langle  x,y\; | \; x^2=y^5=1,
 \; xyx^{-1}=y^{-1} \rangle$ \\
\hline $(4f)$ & $A_4$ & $G(12,3)$ & $(2,3^3)$ & $\langle x,y \; | \;
x=(12)(34), \; y=(123) \rangle$ \\
\hline $(4g)$ & $D_6$ & $G(12,4)$ & $(2^5)$ & $\langle  x,y\; | \; x^2=y^6=1,
 \; xyx^{-1}=y^{-1} \rangle$ \\
\hline $(4h)$ & $D_6$ & $G(12,4)$ & $(2^2,3,6)$ & $\langle  x,y\; | \; x^2=y^6=1,
 \; xyx^{-1}=y^{-1} \rangle$ \\
\hline
$(4i)$ & $D_8$ & $G(16,7)$ & $(2^3, 8)$ & $\langle x,y \; | \; x^2=y^8=1, xyx^{-1}=y^{-1}\rangle$ \\
\hline
$ $ & $ $ & $ $ & $ $ & $ \langle x,y,z,w \; | \; x^2=y^2=z^2=w,$ \\
$(4j)$ & $G(16,9)$ & $G(16,9)$ & $(4^2,8)$ & $w^2=1,\; xzx^{-1}=z^{-1}, $\\
$ $ & $ $ & $ $        & $ $        & $ yzy^{-1}=z^{-1}, \; yxy^{-1}=(xz)^{-1}
\rangle$ \\
\hline $(4k)$ & $\mathbb{Z}_3 \times S_3$ & $G(18,3)$ & $(2^2, 3^2)$ & $\langle z \; | \;
z^3=1 \rangle \times \langle x,y\;|\;x=(12), \;
y=(123) \rangle $ \\
\hline $(4l)$ & $\mathbb{Z}_3 \times S_3$ & $G(18,3)$ & $(3,6^2)$ & $ \langle z \; | \;
z^3=1 \rangle \times \langle x,y\;|\;x=(12), \;
y=(123) \rangle$ \\
\hline
$ $ & $ $ & $ $ & $ $ & $\langle x,y,z \; | \; x^2=y^3=z^3=1,$ \\
$(4m)$ & $\mZ_2 \ltimes (\mZ_3)^2$ & $G(18,4)$ & $(2^2, 3^2)$ &
$xyx^{-1}=y^{-1}, \; xzx^{-1}=z^{-1}, \; [y,z]=1 \rangle $ \\
\hline $(4n)$ & $\mathbb{Z}_2 \times D_5$ & $G(20,4)$ & $(2^3,5)$ & $\langle z \; | \;
z^2=1 \rangle \times \langle  x,y\; | \; x^2=y^5=1,
 \; xyx^{-1}=y^{-1} \rangle$ \\
\hline $(4o)$ & $D_{4,5,-1}$ & $G(20,1)$ & $(4^2,5)$ & $\langle x,y \;|\; x^4=y^5=1,
\; xyx^{-1}=y^{-1}  \rangle$ \\
\hline $(4p)$ & $D_{4,5,2}$ & $G(20,3)$ & $(4^2,5)$ & $\langle x,y \;|\; x^4=y^5=1,
\; xyx^{-1}=y^2  \rangle$ \\
\hline
$(4q)$ & $S_4$ & $G(24,12)$ & $(2^3,4)$ & $\langle x,y \; | \; x=(1234), \; y=(12)   \rangle$ \\
\hline $(4r)$ & $D_{2,12,7}$ & $G(24,10)$ & $(2,6,12)$ & $\langle x,y \;|\; x^2=y^{12}=1,
\; xyx^{-1}=y^7  \rangle$ \\
\hline $(4s)$ & $\textrm{SL}_2(\mathbb{F}_3)$ & $G(24,3)$ & $(3,4,6)$ & $ \langle x,y \;
| \;
 x=
\left(
              \begin{array}{cc}
                1 & 1 \\
                0 & 1 \\
              \end{array}
            \right), \;
y= \left(
              \begin{array}{cc}
                \;\;0 & \;\;1 \\
                -1 & -1 \\
              \end{array}
            \right)
\rangle$ \\
\hline $(4t)$ & $D_{2,16,7}$ & $G(32,19)$ & $(2,4,16)$ & $\langle x,y \;|\; x^2=y^{16}=1,
\; xyx^{-1}=y^7  \rangle$ \\
\hline
 $ $ & $ $ & $ $ & $ $ & $ \langle x,y,z,w \; | \;
x^2=y^2=z^3=w^3=1,
$ \\
$(4u)$ &  $(\mathbb{Z}_2)^2 \ltimes (\mathbb{Z}_3)^2$ & $G(36,10)$ & $(2^3,3)$ & $yzy^{-1}=z^2, \; xwx^{-1}=w^2, $ \\
$ $ & $$ & $ $ & $ $ & $[x,y]=[x,z]=[y,w]=[z,w]=1 \rangle$\\
\hline $ $ & $ $ & $ $ & $ $ & $ \langle x,y,z,w \; | \;
x^2=y^2=z^3=w^3=1,$ \\
$(4v)$ &  $(\mathbb{Z}_2)^2 \ltimes (\mathbb{Z}_3)^2$ & $G(36,10)$ & $(2, 6^2)$ & $yzy^{-1}=z^2, \; xwx^{-1}=w^2, $ \\
$ $ & $$ & $ $ & $ $ & $[x,y]=[x,z]=[y,w]=[z,w]=1 \rangle$\\
\hline
$(4w)$ &  $\mathbb{Z}_6 \times S_3$ & $G(36,12)$ & $(2,6^2)$ & $\langle z~|~z^6=1\rangle\times\langle x,y~|~x=(12),\;y=(123)\rangle$ \\
\hline $(4x)$ &  $\mathbb{Z}_3\times A_4$ & $G(36,11)$ & $(3^2,6)$ & $\langle
z~|~z^3=1\rangle\times\langle x,\,y~|~x=(12)(34),\;y=(123)\rangle$\\
\hline $ $ & $ $ & $ $ & $ $ & $ \langle x,y,z \; | \; x^4=y^3=z^3=1,
$ \\
$(4y)$ &  $\mathbb{Z}_4 \ltimes (\mathbb{Z}_3)^2$ & $G(36,9)$ & $(3,4^2)$ & $xyx^{-1}=yz^2, \; xzx^{-1}=y^2z^2, \; [y,z]=1 \rangle$\\
\hline $ $ & $ $ & $ $ & $ $ & $ \langle x,y,z\; | \; x^2=y^4=z^5=1,
$ \\
$(4z)$ &  $D_4 \ltimes \mZ_5$ & $G(40,8)$ & $(2,4,10)$ & $xyx^{-1}=y^{-1}, \; xzx^{-1}=z, \;
yzy^{-1}=z^{-1} \rangle$ \\
\hline $(4aa)$ & $A_5$ & $G(60,5)$ & $(2, 5^2)$ & $\langle x,y \; | \; x=(12)(34),
\; y=(12345) \rangle$\\
\hline $(4ab)$ &  $\mZ_3\times S_4$ & $G(72,42)$ & $(2,3,12)$ & $\langle
z~|~z^3=1\rangle \times \langle x,\,y~|~x=(12),\;y=(1234)\rangle$ \\
\hline
 $ $ & $ $ & $ $ & $ $ & $ \langle x,y,z,w \; | \;
x^2=y^4=z^3=w^3=1,
$ \\
$(4ac)$ &  $D_4\ltimes(\mZ_3)^2$ & $G(72,40)$ & $(2,4,6)$ &
$xyx^{-1}=y^{-1}, \; xzx^{-1}=w , \; yzy^{-1}=w,$ \\
$ $ & $$ & $ $ & $ $ & $ywy^{-1}=z^2,\;[z,w]=1 \rangle$\\
\hline $(4ad)$ & $S_5$ & $G(120,34)$ & $(2,4,5)$ & $\langle x,y \; | \; x=(12),
\; y=(12345) \rangle$\\
\hline
\end{tabular}
\end{center}
%\textrm{Table \ref{4-nonabelian} (continued on the next page).}
\caption{Nonabelian groups of automorphisms acting with rational
quotient on Riemann surfaces of genus $4$.} \label{4-nonabelian}
\end{table}

\newpage

\section*{Appendix $B$}

The following is the GAP4 script that we used in order
to check that the group $G=\mZ_2 \ltimes (\mZ_2 \times \mZ_4)=G(16,13)$
is not $(1\; | \;2^1)-$generated
(see Proposition \ref{Table-3-no-1-2-gen}, case $(3j)$).
In fact, the output shows that if $[h_1,h_2]$ has order $2$ then
either $\langle h_1, h_2 \rangle \cong G(8,3)=D_4$ or
$\langle h_1, h_2 \rangle \cong G(8,4)=Q_8$.
Completely similar scripts can be used in order to check the other results stated in
 Propositions \ref{no-generation-1-21}, \ref{Table-3-no-1-2-gen},
 \ref{Table-3-no-1-4-gen} and \ref{Table-4-no-1-2-gen}, although in almost all cases
it is also possible to carry out the computations ``by
hand''. \\
\begin{verbatim}
gap> f:=FreeGroup("x", "y", "z");
<free group on the generators [x,y,z]>
gap> x:=f.1; y:=f.2; z:=f.3;
x
y
z
gap> G:=f/[x^2, y^2, z^4,
Comm(x,z), Comm(y,z), x*y*x^-1*(y*z^2)^-1]; # insert the presentation of G
<fp group on the generators [x,y,z]>
gap> x:=G.1; y:=G.2; z:=G.3;
x
y
z
gap> IdSmallGroup(G); # check the IdSmallGroup(G)
[16,13]
gap> for h1 in G do
> for h2 in G do
> H:=Subgroup(G, [h1,h2]);
> if Order(h1*h2*h1^-1*h2^-1)=2 then  # check whether [h1,h2] has order 2
> Print(IdSmallGroup(H), " ");  # identify the subgroup generated by h1 and h2
> fi; od; od; Print("\n");
[8,3] [8,3] [8,3] [8,3] [8,3] [8,3] [8,3] [8,3]
[8,3] [8,3] [8,3] [8,3] [8,3] [8,3] [8,3] [8,3]
[8,3] [8,3] [8,4] [8,3] [8,4] [8,3] [8,4] [8,4]
[8,3] [8,4] [8,4] [8,3] [8,3] [8,4] [8,4] [8,3]
[8,3] [8,3] [8,3] [8,3] [8,3] [8,3] [8,3] [8,3]
[8,3] [8,4] [8,4] [8,3] [8,3] [8,4] [8,4] [8,3]
[8,3] [8,3] [8,3] [8,3] [8,3] [8,3] [8,3] [8,3]
[8,3] [8,3] [8,3] [8,3] [8,3] [8,3] [8,3] [8,3]
[8,3] [8,3] [8,4] [8,3] [8,4] [8,3] [8,4] [8,4]
[8,3] [8,4] [8,4] [8,3] [8,3] [8,4] [8,4] [8,3]
[8,3] [8,4] [8,4] [8,3] [8,3] [8,4] [8,4] [8,3]
[8,3] [8,3] [8,3] [8,3] [8,3] [8,3] [8,3] [8,3]
gap>
\end{verbatim}

\end{document}